%
%

\input amstex  

\input amssym
\input amssym.def

\magnification 1200
\loadmsbm
\parindent 0 cm


\define\nl{\bigskip\item{}}
\define\snl{\smallskip\item{}}
\define\inspr #1{\parindent=20pt\bigskip\bf\item{#1}}
\define\iinspr #1{\parindent=27pt\bigskip\bf\item{#1}}
\define\ainspr #1{\parindent=24pt\bigskip\bf\item{#1}}
\define\aiinspr #1{\parindent=31pt\bigskip\bf\item{#1}}

\define\einspr{\parindent=0cm\bigskip}

\define\ot{\otimes}



\centerline{\bf Weak Multiplier Hopf Algebras}
\snl\centerline{Preliminaries, motivation and basic examples}
\bigskip\bigskip
\centerline{\it  Alfons Van Daele \rm $^{(1)}$ and \it Shuanhong Wang \rm $^{(2)}$}
\bigskip\bigskip\bigskip
{\bf Abstract} 
\nl 
Let $G$ be a {\it finite group}. Consider the algebra $A$ of all complex functions on $G$ (with pointwise product). Define a coproduct $\Delta$ on $A$ by $\Delta(f)(p,q)=f(pq)$ where $f\in A$ and $p,q\in G$. Then $(A,\Delta)$ is a Hopf algebra. If $G$ is only a {\it groupoid}, so that the product of two elements is not always defined, one still can consider $A$ and define $\Delta(f)(p,q)$ as above when $pq$ is defined. If we let $\Delta(f)(p,q)=0$ otherwise, we still get a coproduct on $A$, but $\Delta(1)$ will no longer be the identity in $A\ot A$. The pair $(A,\Delta)$ is not a Hopf algebra but a weak Hopf algebra. If $G$ is a {\it group}, but {\it no longer finite}, one takes for $A$ the algebra of functions with finite support. Then $A$ has no identity and $(A,\Delta)$ is not a Hopf algebra but a multiplier Hopf algebra. Finally, if $G$ is a {\it groupoid}, but {\it not necessarily finite}, the standard construction above, will give, what we call in this paper, a weak multiplier Hopf algebra. 
\snl
Indeed, this paper is devoted to the development of this 'missing link': {\it weak multiplier Hopf algebras}. We spend a great part of this paper to the motivation of our notion and to explain where the various assumptions come from. The goal is to obtain a good definition of a weak multiplier Hopf algebra. Throughout the paper, we consider the basic examples  and use them, as far as this is possible, to illustrate what we do. In particular, we think of the
 finite-dimensional weak Hopf algebras. On the other hand however, we are also inspired by the far more complicated existing analytical theory.
\snl
In our forthcoming papers on the subject, we develop the theory further. In [VD-W2] we start from the definition of a weak multiplier Hopf algebra as it is obtained and motivated in this paper and we prove the main properties. In [VD-W3] we continu with the study of the source and target algebras and the corresponding source and target maps. In that paper, we also give more examples. Finally, in [VD-W4] we study integrals on weak multiplier Hopf algebras and duality. Other aspects of the theory will be considered later.
\nl \nl

{\it October 2012} (Version 1.5)
\vskip 1 cm
\hrule
\bigskip
\parindent 0.7 cm
\item{($1$)} Department of Mathematics, University of Leuven, Celestijnenlaan 200B,
B-3001 Heverlee, Belgium. {\it E-mail}: Alfons.VanDaele\@wis.kuleuven.be
\item{($2$)} Department of Mathematics, Southeast University, Nanjing 210096, China. {\it E-mail}:  Shuanhwang2002\@yahoo.com {\it or} Shuanhwang\@seu.edu.cn
\parindent 0 cm

\newpage



\bf 0. Introduction \rm
\nl
Let $G$ be a {\it groupoid}. It is a set with a distinguished subset of pairs $(p,q)$ in $G\times G$ for which the product $pq$ in $G$ is defined. This product is  associative, in the appropriate sense. The product $pq$ is only defined when the so-called source $s(p)$ of the element $p$ is equal to the target $t(q)$ of the element $q$. The {\it source} and {\it target maps} are defined from $G$ to the set of units and this set can (and will) be considered as a subset of $G$. 
\snl
For a more precise definition and the elementary theory of groupoids, we refer to [Br], [H] and [R].
\nl
\it Examples of groupoids \rm
\nl
Perhaps the simplest example of a groupoid is obtained from an equivalence relation $\sim$ on a set $X$. The elements of the groupoid $G$ are pairs $(y,x)$ with $x,y\in X$ and $x\sim y$. The set of units is $X$ and the source and target maps are given by
$$s(y,x)=x \qquad \text{and} \qquad t(y,x)=y.$$
for $(y,x)$ in $G$. The set of units is considered as a subset of $G$ via the map $x\mapsto (x,x)$. The product of  $(z,y)$ with $(y,x)$ is $(z,x)$ when $x,y,z\in X$ and $x\sim y$ and $y\sim z$.
\snl
This example is of course rather trivial and it will be of little use for illustrating our theory. Instead, the following, also well-known but more involved example is more interesting.
\snl
Now $X$ is a set and $H$ is a group acting on the set $X$, say from the left. We denote by $hx$ the element obtained by letting $h\in H$ act on $x\in X$. We let
$$G=\{(y,h,x)\mid x,y\in X, h\in H, y=hx\}.$$
The space of units is again $X$ and the source and target maps are given by 
$$s(y,h,x)=x \qquad \text{and} \qquad t(y,h,x)=y.$$
The space of units is a subset of $G$ via the map $x\mapsto (x,e,x)$ where $e$ is the identity element in the group $H$. The product of two elements $(z,k,y)$ and $(y,h,x)$ is $(z,kh,x)$ where of course  $kh$ is the product  of the elements $k,h\in H$ in the group $H$. 
\snl
This second example reduces to a special case of the first one if the group $H$ is trivial. On the other hand, the action of $H$ defines an equivalence relation on $X$ and this yields a map from the groupoid in the second example to the groupoid in the first example. Finally, if the set $X$ consists of a single point, the second example will be the same as the group $H$ itself.
\snl
It is good to have these various cases in mind. However, mostly we will only use the case of a general groupoid to illustrate definitions and results in our paper here.
\nl
\newpage 

\it The algebra $K(G)$ and the coproduct on this algebra \rm
\nl
In what follows, let $G$ be any groupoid. Consider the algebra $K(G)$ of complex functions with finite support (and pointwise operations). We will denote this algebra by $A$. Observe that it only has an identity when $G$ is finite but that the product in $A$ is always non-degenerate. The product in $G$ yields a coproduct $\Delta$ on $A$ by the formula
$$
\Delta(f)(p,q)=
\cases
f(pq) & \text{if $pq$ is defined},\\
0 & \text{otherwise}.
\endcases
$$
It is a homomorphism from $A$ to the multiplier algebra $M(A\ot A)$ of the tensor product $A\ot A$ of $A$ with itself. Recall that in this case we have a natural identification of $A\ot A$ with $K(G\times G)$ and of $M(A\ot A)$ with $C(G\times G)$, the algebra of all complex functions on $G\times G$. The associativity of the product will yield the coassociativity of the coproduct (to be understood in the appropriate sense, see Definition 1.1 in Section 1). 
\nl 
If $G$ is {\it finite}, then $A$ has an identity, we have $M(A\ot A)=A\ot A$ and so $\Delta$ will map $A$ to $A\ot A$. There is no difficulty with the notion of coassociativity in this case. However, in general, we will not have that $\Delta(1)=1\ot 1$, precisely because the product in $G$ is not defined everywhere. Instead, $\Delta(1)$ will be an idempotent $E$ in $A\ot A$ satisfying, among other properties, 
$$E\Delta(a)=\Delta(a)
	\qquad \qquad \text{and} \qquad\qquad
		\Delta(a)E=\Delta(a) \tag"(0.1)"$$
for all $a$ in $A$.
In this case, the pair $(A,\Delta)$ is  a (finite-dimensional) weak Hopf algebra in the sense of [B-N-S]. This example is also considered further at various places in this paper.
\snl
If $G$ is no longer assumed to be finite, one expects that the resulting pair $(A,\Delta)$ will be a {\it weak multiplier Hopf algebra}. A weak multiplier Hopf algebra should relate to a (finite-dimensional) weak Hopf algebra as a multiplier Hopf algebra to a (finite-dimensional) Hopf algebra. The new notion should generalize the concept of a (finite-dimensional) weak Hopf algebra to the case of possibly non-unital, infinite-dimensional algebras.
\nl
\it Development of the notion of a weak multiplier Hopf algebra \rm 
\nl
This is exactly what this paper is about. We develop the notion of a weak multiplier Hopf algebra. And we do it in {\it a special way}. We will not start with stating the correct definition with its various objects and conditions and continue with proving properties. Instead, we will gradually develop the notion and spend a great deal of this paper to the motivation of the concept and explain where the various conditions come from. By doing so, we end up with a notion and conditions that will be experienced as very natural. Moreover, this will provide the reader with a better understanding of what weak multiplier Hopf algebras really are. In fact, we also believe that this work might contribute even to a better understanding of the existing theory of weak Hopf algebras as well as provide guidelines to find a possibly simpler approach to the present (rather complicated) treatment of locally compact quantum groupoids and their measured counterparts.
\snl
In a second paper on the subject [VD-W2], we will follow the traditional path and start with the definition as we find it here and develop the theory from this definition. In a third paper [VD-W3] we will continue this study and obtain more properties of the source and target algebras, as well as of the source and target maps. In that paper, we will also give more examples. In a fourth paper [VD-W4], we treat integrals on weak multiplier Hopf algebras and duality.
\snl
Our approach to weak multiplier Hopf algebras in this paper is inspired by the theory of multiplier Hopf algebras as developed in the original paper on this theory [VD1]. We will say more about this in the forthcoming item in this introduction where we give an overview of the content of the paper.
\snl
We use the motivating examples, coming from a groupoid, to illustrate various aspects of our approach. Of course, another main source of inspiration is the existing theory of weak Hopf algebras (as developed by G.\ B\"ohm, F.\ Nill \& K.\ Szlach\'anyi in [B-N-S] and [B-S]). 
\snl
All of these cases will help to obtain the right axioms and conditions for a weak multiplier Hopf algebra. There is however still another criterion we will use to see if the conditions are good. They should be as much as possible 'self-dual'. By this we mean that, for a dual pair, a given condition on one component should yield a given condition on the other component by duality. The idea will be more clear when we come to the first application of this principle (see e.g.\ Remark 2.3 in Section 2).
\nl
\it The idempotent $E$ playing the role of $\Delta(1)$ \rm
\nl
Before we start with the description of the content of the paper, let us discuss what is {\it the first difficulty} we encounter as a consequence of the fact that we do not assume our algebra to be unital.
\snl
If we look at the definition of a weak Hopf algebra $(A,\Delta)$, we see that several conditions are formulated in terms of the idempotent $\Delta(1)$ in $A\ot A$. In the case we consider here, we simply can not consider $\Delta(1)$ in the first place.
\snl
If we take a multiplier Hopf algebra $(A,\Delta)$, there is a possibility to extend the coproduct $\Delta$ to a homomorphism from $M(A)$ to $M(A\ot A)$ using the fact that the coproduct is assumed to be non-degenerate. This means that $\Delta(1)$ is defined but on the other hand, it is necessarily equal to $1\ot 1$ in $M(A\ot A)$. 
\snl
In the case of weak multiplier Hopf algebras, we can not expect this to happen and so, in particular, we can not assume that the coproduct is non-degenerate. Hence we can not apply this extension procedure and obtain $\Delta(1)$ in this way.
\snl
Fortunately, there is a weaker notion of non-degeneracy that can be used in this new context. It does allow to extend $\Delta$ to $M(A)$ and hence give a meaning to $\Delta(1)$. It will be an idempotent $E$ in $M(A\ot A)$ satisfying the formulas (0.1) above. By the conditions, it will be uniquely defined. This is explained in an appendix of this paper. See also further in this introduction when we describe the content of the paper.
\snl
We believe that the solution we have here for this difficulty is precisely what makes a satisfactory theory of weak multiplier Hopf algebras possible.
\nl
\it The content of the paper \rm
\nl
In the {\it first section}, we start with an algebra $A$ and we give the precise definition of a coproduct $\Delta$ as we will use it in this theory. We recall the definition of a counit and show that this is only a good notion if we also assume that the coproduct is full, meaning roughly speaking that the legs of $\Delta$ are all of $A$. We have to be more careful because, as mentioned already earlier,  we can not expect the coproduct to be non-degenerate as is usually done in the case of multiplier Hopf algebras.  Another problem with the counit comes from the fact that it is no longer expected to be a homomorphism.
\snl
We consider also coproducts on $^*$-algebras, as well as regular coproducts in general. 
\snl
We give a preliminary definition of a dual pair in this context. We will not really investigate this notion, but rather use it to motivate the conditions we impose throughout the development further in this paper. As mentioned already, we want these conditions to be 'self-dual'.
\snl
In the first section, we also give {\it elementary examples}.
\snl
In the {\it second section} we treat the antipode(s). The starting point is a pair of generalized inverses $R_1$ and $R_2$ of the canonical maps $T_1$ and $T_2$, defined from $A\ot A$ to itself by 
$$T_1(a\ot b)=\Delta(a)(1\ot b)\qquad\qquad\text{and}\qquad\qquad T_2(a\ot b)=(a\ot 1)\Delta(b).$$
Under some natural conditions, the inverses yield antipodes $S_1$ and $S_2$ such that $R_1$ and $R_2$ are given by the well-known formulas in terms of these antipodes. The conditions are 'self-dual' in the sense mentioned earlier. In the regular case, we have two other canonical maps $T_3$ and $T_4$ and a good choice of generalized inverses $R_3$ and $R_4$ yield two other (inverse) antipodes $S_3$ and $S_4$. In the $^*$-case, regularity of the coproduct is automatic and we find natural relations between the pair $(S_3,S_4)$ and the pair $(S_1,S_2)$. 
\snl
The rest of the section is devoted to all sorts of properties of these two (respectively four) antipodal maps and relations among them. Most of these further results are included mainly for motivational reasons.
\snl
At the end of the section, we formulate the notion of a {\it unifying multiplier Hopf algebra}. It is by no means yet the final object (weak multiplier Hopf algebras) and for the moment, we do not plan to study the concept further.
\snl
We illustrate all this using examples.
\snl
In {\it Section} 3, we will impose some extra conditions. First, we look at the ranges of the maps $T_1$ and $T_2$. Essentially we will require that there is an idempotent $E$ in $M(A\ot A)$ such that
$$E(A\ot A)=\Delta(A)(1\ot A) 
	\quad \qquad \text{and} \quad\qquad
         (A\ot A)E=(A\ot 1)\Delta(A),
$$
together with the possibility of making choices for $R_1$ and $R_2$ that take into account this property. As we  show in an appendix, this will imply that we are able to extend the map $\Delta$, as well as the maps $\iota\ot\Delta$ and $\Delta\ot\iota$, to homomorphisms on $M(A)$, respectively $M(A\ot A)$ (where $\iota$ is used for the identity map). We will have that $\Delta(1)=E$ as expected. 
\snl 
Further in this section, we consider the kernels of the maps $T_1$ and $T_2$. Also here, we impose some natural conditions so that we have some equivalent of the multiplier $E$ related with these kernels instead of the ranges. We investigate relations between these objects and find that the equality $S_1=S_2$ of the antipodes found in Section 2, is very basic. In fact, we see that this requirement fundamentally fixes these antipodes (and so also the kernels of the maps $T_1$ and $T_2$), given the multiplier $E$ above. We feel that this is quite remarkable (see a remark in Section 3).
\snl
In {\it Section} 4, we finally come up with a {\it first definition} of a weak multiplier Hopf algebra, based on the considerations of the previous sections. It is a definition based on the characterizations of the kernels and the ranges of the canonical maps $T_1$ and $T_2$. The definition does not involve the antipode, just like in the original definition of a multiplier Hopf algebra. It is shown that the basic examples, coming from a groupoid, satisfy the axioms. Also a weak Hopf algebra is a weak multiplier Hopf algebra with our definition.
\snl
Further in Section 4, we discuss the {\it regular case}. A weak multiplier Hopf algebra is defined regular if the antipode is a bijective map from the algebra $A$ to itself. We also briefly treat the involutive case here. In this section, some of results are not proven and we refer to [VD-W2] for details.
\snl
Finally, in the last section, {\it Section 5}, we formulate some conclusions and indicate what will be done in the next papers on the subject (see also further). 
\nl
As mentioned already, in {\it Appendix} A we study the problem of extending the coproduct $\Delta$, as well as the derived maps $\Delta\ot \iota$ and $\iota\ot\Delta$ to the multiplier algebras $M(A)$ and $M(A\ot A)$ respectively. Remember that in this theory, we can not require that $\Delta$ is non-degenerate and so we can not use the existing method. Instead, we will use the existence of the smallest idempotent multiplier $E$ satisfying $E\Delta(a)=\Delta(a)$ and $\Delta(a)E=\Delta(a)$ for all $a$. As we will see in Section 3,  this is something we can assume instead of non-degeneracy and it still works to get nice extensions.
\nl
The further development of the theory is continued in separate papers on the subject. In the first one [VD-W2], we start from the definition of a weak multiplier Hopf algebra as we have found it here and develop the theory from the definition. In a way, this is first reorganizing the material from this paper, leaving out the motivational part and add some proofs that are omitted in the treatment here. 
In the second paper [VD-W3] we investigate further properties of the source and target maps and of the source and target algebras. We then are ready for the construction of some more examples. Finally, in [VD-W4], we will study integrals on weak multiplier Hopf algebras and the related theory of duality.
\nl
\it Conventions and notations \rm
\nl
All our algebras will be algebras over the field $\Bbb C$ of complex numbers (although probably algebras over other fields will be possible as well). We do not assume that they have a unit. But we do assume that the product, seen as a bilinear map, is non-degenerate. When $A$ is such an algebra, we will denote by $M(A)$ the multiplier algebra of $A$. It is characterized as the largest algebra with $1$ containing $A$ as a dense ideal. 
\snl
In general we also will need the algebras to be idempotent in the sense that any element is a sum of products of elements in the algebra. We write this as $A^2=A$. In fact, we will see later (see Proposition 4.9 in [VD-W2]) that the underlying algebra of a regular weak multiplier Hopf algebra automatically has local units. This of course implies that the product is non-degenerate and that the algebra is idempotent.
\snl
Occasionally, we will need to work with the algebra $L(A)$ of left multipliers of $A$ and the algebra $R(A)$ of right multipliers of $A$. Elements of $L(A)$ are maps $\lambda:A\to A$ satisfying $\lambda(ab)=\lambda(a)b$ for all $a,b\in A$. Similarly elements of $R(A)$ are maps $\rho:A\to A$ satisfying $\rho(ab)=a\rho(b)$ for all $a,b\in A$. When $m\in L(A)$ we will write $ma$ instead of $m(a)$ for $a\in A$. When on the other hand $m\in R(A)$, we write $am$ for $m(a)$. This notation is compatible with the defining properties and consistent in the case where $M(A)$ is considered as a subset of both $L(A)$ and $R(A)$. Remark that $M(A)$, $L(A)$ and $R(A)$ are all identified with $A$ itself in the case of a unital algebra. 
\snl
We use $1$ to denote the identity element in the algebra $M(A)$. We use $\iota$ for the identity map on any vector space. Finally we use $e$ for the identity element in a group.
\snl
When $T$ is a linear map from one vector space to another, we will denote the kernel of $T$ by $\text{Ker}(T)$ and the range of $T$ by $\text{Ran}(T)$.
\snl
We will use $A^{\text{op}}$ for the algebra $A$ with the opposite product. And when $\Delta$ is a coproduct on $A$, we use $\Delta^{\text{cop}}$ for the new coproduct on $A$ obtained by composing $\Delta$ with the flip map (extended to the multiplier algebra).
\snl
We will use the {\it leg-numbering} notation in different circumstances. If e.g.\ $A$ is an algebra (with or without identity) and $x\in A\ot A$, we will use $x_{12}$, $x_{13}$ and $x_{23}$ for the elements $x$ as sitting in the multiplier algebra of the threefold tensor product $A\ot A\ot A$ in the three possible ways. We have $x_{12}=x\ot 1$, $x_{23}=1\ot x$ and $x_{13}=\sigma_{23}(x\ot 1)$ where $\sigma_{23}$ flips the last two factors in the threefold tensor product. A similar convention is used for maps with the identity element of the algebra replaced by the identity map $\iota$. This notation is also used for coproducts in the obvious sense.
\snl
We will often use the Sweedler notation for a coproduct $\Delta$. In the case of a coproduct in the usual sense, so a map from a vector space $A$  to the tensor product $A\ot A$, this presents no problem. As in our case, the coproduct mostly does not map into this tensor product, but rather into a bigger space, the Sweedler notation has to be used with care. Sufficient {\it covering} is needed. For a detailed account on the use of the Sweedler notation in this more general context, we refer to [VD3].
\snl
Finally we use the items {\it Definition}, {\it Proposition}, etc. in the usual sense. But in this paper occasionally, we will also have items like {\it Condition} and {\it Assumption}. In the first case, it will be a condition that we temporarily assume (and we will mention when we do this) whereas in the second case, we formulate a condition that will be assumed from then onwards. We will recall this convention when we first use these terms.

\newpage

\it Basic references \rm
\nl
For the theory of Hopf algebras, we refer to the basic works of Abe [A] and Sweedler [S]. For multiplier Hopf algebra the references are [VD1] and [VD2]. For the use of the Sweedler notation for coproducts with values in the multiplier algebras and the covering technique, a good and detailed account is found in [VD3]. We have a note that might be useful when dealing with local units and multiplier algebras ([VD-Ve]).
\snl
The theory of weak Hopf algebras is developed in [B-N-S] and [B-S] by G.\ B\"ohm, F.\ Nill \& K.\ Szlach\'anyi. However, our notations will be more inspired by the survey paper by D.\ Nikshych \& L.\ Vainerman [N-V2]. We also refer to various papers in the proceedings of a conference in Strasbourg on locally compact quantum groups and groupoids [V]. In these proceedings, also material is found about the analytical theory of locally compact and measured quantum groupoids.
\nl
\bf Acknowledgements \rm
\nl
The first named author (Alfons Van Daele) would like to express his thanks to the following people. He is greatly indebted to Leonid Vainerman for introducing him to the subject of weak Hopf algebras and to his coauthor Shuanhong Wang for motivating him to start the research on weak multiplier Hopf algebras. He is also grateful to Michel Enock and Michel Vallin for informal discussions on the topological theory of locally compact quantum groupoids (or measured quantum groupoids).
\snl
The second named author (Shuanhong Wang) would like to thank his coauthor,
for his help and advices when visiting the Department of Mathematics of the K.U. Leuven in Belgium during several periods in 2004, 2006 and 2009 as a post-doctoral fellowship. He also is grateful to the analysis research group of the K.U. Leuven for the warm hospitality. In particular, he thanks Johan Quaegebeur, Stefaan Vaes and Johan Kustermans for all their help on his life and research during his 
stay in Leuven. This work is part of a project, supported by a research fellowship from the K.U. Leuven in 2004. 
\nl\nl



\bf 1. The coproduct and the counit \rm
\nl
Let $A$ be an (associative) algebra over the field $\Bbb C$ of the complex numbers. We do not assume that $A$ has a unit, but we do require that the product, seen as a bilinear form, is non-degenerate. This means that, whenever $a\in A$ and $ab=0$ for all $b\in A$ or $ba=0$ for all $b\in A$, we must have that $a=0$. Then we can consider the multiplier algebra $M(A)$ of $A$. Recall that $M(A)$ is the largest algebra with identity containing $A$ as a dense two-sided ideal. In particular, we still have that,
whenever $a\in M(A)$ and $ab=0$ for all $b\in A$ or $ba=0$ for all $b\in A$, again $a=0$. Further we consider the algebra $A\ot A$. It is still non-degenerate and we also have its multiplier algebra $M(A\ot A)$. There are natural embeddings 
$$A\ot A \subseteq M(A)\ot M(A) \subseteq M(A\ot A).$$
Observe that in general, when $A$ has no identity, these two inclusions are strict.
\snl
Later, we will need to assume that our algebra $A$ is also idempotent, that is that any element in $A$ is a sum of products of elements in $A$. We write this condition as $A=A^2$. There are reasons to believe that this condition follows from the others, but we have not been able to prove this. For regular weak multiplier Hopf algebras, it will be shown in [VD-W2] that local units exist. This is of course a stronger property than being idempotent. There are also reasons to believe that the existence of local units is also true in the non-regular case, but again we have not been able to show this.
\snl
We will also consider $^*$-algebras. Then $A$ has a conjugate linear involution $a\mapsto a^*$ satisfying $(ab)^*=b^*a^*$ for all $a,b\in A$. In this case, also $M(A)$ is a $^*$-algebra in a natural way. The same is true for $A\ot A$ and $M(A\ot A)$.
\snl
We have the following (basic) example in mind. Let $X$ be a set and let $A$ be the algebra $K(X)$ of complex functions on $X$ with finite support. For the pointwise product, it is a non-degenerate algebra. Clearly $A\ot A$ is naturally identified with $K(X\times X)$ while $M(A)$ and $M(A\ot A)$ are identified with $C(X)$ and with $C(X\times X)$, the algebras of all complex functions on $X$  and on $X\times X$ respectively. We have a natural $^*$-structure on this algebra. It is given by $f^*(x)=\overline{f(x)}$ whenever $f\in K(X)$ and $x\in X$.
%
\nl
\it The coproduct \rm
\nl
Now we are ready to introduce the concept of a comultiplication (or coproduct) as we will use it in this theory.

\inspr{1.1} Definition \rm
Let $A$ be an algebra with a non-degenerate product. A {\it coproduct} on $A$ is a homomorphism $\Delta:A\to M(A\ot A)$ satisfying 
\snl
i) $\Delta(a)(1\ot b)$ and $(a\ot 1)\Delta(b)$ are in $A\ot A$ for all $a,b\in A$, \newline
ii) $\Delta$ is {\it coassociative} in the sense that
$$(a\ot 1\ot 1)(\Delta \ot \iota)(\Delta(b)(1\ot c))=
      (\iota \ot \Delta)((a\ot 1)\Delta(b))(1\ot 1\ot c)$$
for all $a,b,c\in A$, where $\iota$ is used for the identity map on $A$. In the case of a $^*$-algebra we assume that $\Delta$ is a $^*$-homomorphism.
\hfill$\square$\einspr

Remark that i) is needed to formulate ii). This is also the way coassociativity was first introduced for multiplier Hopf algebras (see [VD1]). Other ways to formulate coassociativity, used for multiplier Hopf algebras are based on the non-degeneracy of the coproduct, making it possible to extend the homomorphisms $\Delta\ot \iota$ and $\iota\ot \Delta$ to $M(A)$. Then coassociativity can simply be written in its usual form, i.e.\ as $(\Delta\ot\iota)\Delta=(\iota\ot\Delta)\Delta$. This however is not so simple here as we will not have that $\Delta$ is non-degenerate. There is a workaround, but we need more conditions. We discussed this already in the introduction and we will come back to this problem in Section 3 (as well as in the appendix). 
\snl
On the other hand, the conditions i) are very natural as we know from the theory of multiplier Hopf algebras and as we will see shortly when we consider our basic examples. The conditions allow us to define the following {\it canonical maps}, playing a very important role in our approach to the theory.

\inspr{1.2} Notation \rm
When $\Delta$ is a coproduct on $A$ as in the previous definition, we denote by $T_1$ and $T_2$, the linear maps from $A\ot A$ to itself, given by 
$$T_1(a\ot b)=\Delta(a)(1\ot b)\quad\qquad\text{and}\quad\qquad T_2(a\ot b)=(a\ot 1)\Delta(b)$$
for $a,b\in A$. 
\hfill$\square$\einspr

Remark that coassociativity, as defined in Definition 1.1, can be written with the help of these maps as 
$$(T_2\ot \iota)(\iota\ot T_1)=(\iota\ot T_1)(T_2\ot \iota).$$
\snl
In general, we will not automatically have that also the elements $\Delta(a)(b\ot 1)$ and $(1\ot a)\Delta(b)$, defined in the multiplier algebra $M(A\ot A)$, will belong to $A\ot A$. This takes us to the following definition and notation.

\inspr{1.3} Definition \rm
A coproduct $\Delta$ on a non-degenerate algebra $A$ is called {\it regular} if also 
\snl
iii) $\Delta(a)(b\ot 1)$ and $(1\ot a)\Delta(b)$ belong to $A\ot A$ for all $a,b\in A$.
\hfill$\square$\einspr

If $\Delta$ is regular, then $\Delta^{\text{cop}}$, obtained from $\Delta$ by composing it with the flip map, is a coproduct in the sense of Definition 1.1. Also remark that regularity is automatic in the case of a coproduct on a $^*$-algebra because $\Delta$ is then assumed to be a $^*$-homomorphism. It is also automatic if the algebra is abelian or if the coproduct is coabelian (i.e\ when $\Delta^{\text{cop}}=\Delta$).
\snl
In the case of a regular coproduct, we also introduce the following notations.

\inspr{1.4} Notation \rm
For a regular coproduct $\Delta$ on a non-degenerate algebra $A$, we denote by $T_3$ and $T_4$ the linear maps from $A\ot A$ to itself, given by
$$T_3(a\ot b)=(1\ot b)\Delta(a)\qquad\quad \text{and} \qquad\quad T_4(a\ot b)=\Delta(b)(a\ot 1)$$
for $a,b\in A$.
\hfill$\square$\einspr

The notational conventions are such that, in the case of a $^*$-algebra, the pair $(T_1,T_2)$ is converted to the pair $(T_3,T_4)$ by the involution. Other conventions are possible, e.g.\ the ones that are suggested by duality (cf.\ the following item). 
\nl
\it Dual pairs \rm
\nl
This is the point where we can give a {\it preliminary definition of a dual pair}. A more complete definition will be given in [VD-W4], but for our purposes, this is not yet needed.

\inspr{1.5} Definition \rm
Let $(A,\Delta)$ and $(B,\Delta)$ be two pairs of non-degenerate algebras with a coproduct as in Definition 1.1. We use the same symbols for the maps defined in Notation 1.2 for each of them. A pairing between $(A,\Delta)$ and $(B,\Delta)$ is a non-degenerate bilinear map $(a,b)\mapsto \langle a,b\rangle$ from $A\times B$ to $\Bbb C$ such that
$$\align \langle T_1(a\ot a'),b\ot b'\rangle &= \langle a\ot a', T_2(b\ot b') \rangle \\
			\langle T_2(a\ot a'),b\ot b'\rangle &= \langle a\ot a', T_1(b\ot b') \rangle
\endalign$$
for all $a,a'\in A$ and $b,b'\in B$. 
\hfill$\square$\einspr

If the algebras have identities, we see that these two conditions are equivalent with the fact that the coproduct on one algebra is dual to the product on the other one.
\snl
If the two pairs are regular, we also will have
$$\align \langle T_3(a\ot a'),b\ot b'\rangle &= \langle a'\ot a, T_3(b'\ot b) \rangle \\
			\langle T_4(a\ot a'),b\ot b'\rangle &= \langle a'\ot a, T_4(b'\ot b) \rangle
\endalign$$
for all $a,a'\in A$ and $b,b'\in B$. Observe that these equations involve flip maps and do not convert $T_3$ to $T_4$ as in the case of $T_1$ and $T_2$. This is due to the choice we made in Notation 1.4 where we have given priority to the implications of having an involution. Another choice would have made these formulas nicer.
\snl
At this moment, we can not yet give the correct definition of a pairing in the case of $^*$-algebras as we will need the antipode to do this.
\snl
It is this notion of duality that we will be using throughout the paper {\it mainly for motivational reasons}. As mentioned already, in a forthcoming paper on the subject, we will develop, fine tune and study the notion in greater detail (see [VD-W4]).
\nl
Up to now, nothing prevents $\Delta$ to be completely trivial (i.e.\ identically $0$). Of course, we will exclude this case.
\snl
But first, we consider the {\it examples associated with a groupoid}. We have the following two propositions.

\inspr{1.6} Proposition \rm
Let $G$ be a groupoid. Let $A$ be the algebra $K(G)$ of complex functions with finite support on $G$, with pointwise product. Let $\Delta$ be defined by
$$
\Delta(f)(p,q)=
\cases
f(pq) & \text{if $pq$ is defined},\\
0 & \text{otherwise}.
\endcases
$$
Then $\Delta$ is a regular coproduct on $A$. If moreover $A$ is considered with its natural $^*$-structure, given by $f^*(p)=\overline{f(p)}$ when $f\in K(G)$ and $p\in G$, then $\Delta$ is a $^*$-homomorphism.

\snl\bf Proof\rm:
It is clear that $\Delta$ is a homomorphism from $A$ to $M(A\ot A)$ and a $^*$-homo\-mor\-phism for the natural $^*$-operation defined on $A$.
\snl
Take $f,g\in A$ and consider the function $\Delta(f)(1\ot g)$. It maps the pair $(p,q)$ to $f(pq)g(q)$ if $pq$ is defined and to $0$ otherwise. The presence of $g$ forces $q$ to lie in a finite set (for the result to be non-zero). Also $pq$ must lie in a finite set and because $p=(pq)q^{-1}$ when $pq$ is defined, the result will be $0$ except when also $p$ lies in a finite set. Therefore $\Delta(f)(1\ot g)\in K(G\times G)$. Similarly for $(f\ot 1)\Delta(g)$ so that condition i) in Definition 1.1 is satisfied. 
\snl
The coassociativity of $\Delta$, as formulated in ii) of the definition, is a straightforward consequence of the associativity of the product in $G$. Regularity here follows automatically because the algebra is abelian.
\hfill $\square$
\einspr 

The identity for this algebra is the function that sends all elements of $G$ to $1$. This will only belong to $A$ if $G$ is finite.  We have however a {\it natural candidate} for $\Delta(1)$. It should be the function on $G\times G$ that maps a pair $(p,q)$ to $1$ if $pq$ is defined and to $0$ otherwise. This is a self-adjoint idempotent in $C(G\times G)$. And of course 
$\Delta(1)\Delta(f)=\Delta(f)$ for all $f\in K(G)$. It is the smallest idempotent with this property. In Section 3, we will consider such an idempotent in general.
\snl
Observe also that $A^2=A$ for this algebra.
\nl
In the following proposition, we consider {\it the dual case}.

\inspr{1.7} Proposition \rm
Let $G$ be a groupoid. Let $B$ be the space $\Bbb C G$ of complex functions with finite support on $G$.  It is a non-degenerate associative algebra for the convolution product (defined below). We use $p\mapsto\lambda_p$ for the imbedding of $G$ in $\Bbb C G$. If then we define $\Delta$ on $B$ by $\Delta(\lambda_p)=\lambda_p\ot \lambda_p$ for all $p\in G$, we get a regular coproduct on $B$. If moreover $B$ is considered with its natural $^*$-structure, given by $\lambda_p^*=\lambda_{p^{-1}}$, when $p\in G$, then $\Delta$ is a $^*$-homomorphism.

\snl\bf Proof\rm: The convolution product is defined by $\lambda_p\lambda_q=\lambda_{pq}$ when $pq$ is defined and $\lambda_p\lambda_q=0$ otherwise. It makes $B$ into an associative algebra (because of the associativity of the product in $G$). It need not have an identity (see a remark following this proof) but the product is non-degenerate. To see this, let $a=\sum a(p)\lambda_p$ be an element in $B$ and assume that $ab=0$ for all $b$. Fix  $p_0$, let $q$ be the source of $p_0$ in $G$ and put $b=\lambda_q$. Then $\sum a(p)\lambda_{pq}=0$ where now the sum is only taken over those elements $p$ in $G$ that have the same source as $p_0$. This will imply that $a(p)=0$ for all these elements, in particular $a(p_0)=0$. Therefore $a=0$.
\snl
In this case, we have $\Delta(B)\subseteq B\ot B$ and the conditions i) and iii) in the Definitions 1.1 and 1.3 are automatic. Also coassociativity is straightforward. The same is true for the last statement about the involutive structure.
\hfill $\square$
\einspr 

The identity in the multiplier algebra is $\sum \lambda_e$ where the sum is taken over all units $e$ in $G$. This will belong to $B$ if and only if the set of units is finite. Otherwise, it will not and then we get an algebra without identity.
\snl
Also here, there is a natural candidate for $\Delta(1)$, namely $\sum \lambda_e\ot \lambda_e$ where again the sum is taken over the set of units in $G$. It satisfies $\Delta(1)\Delta(\lambda_p)=\Delta(\lambda_p)$ and $\Delta(\lambda_p)\Delta(1)=\Delta(\lambda_p)$ for all $p\in G$ and also here, it is the smallest idempotent in $M(B\ot B)$ with this property.  
\snl
The algebra $B$ is also idempotent.
\snl
The two cases, given in Proposition 1.6 and 1.7 are really dual to each other in the sense of Definition 1.5. The pairing is the obvious one.
\snl
At the end of this section, we will consider these examples again and also the one coming from a weak Hopf algebra.
\nl
\it Full coproducts and the counit \rm
\nl
First we consider the concept of a counit. The following seems to be a natural definition.

\inspr{1.8} Definition \rm
Let $A$ be a non-degenerate algebra and $\Delta$ a coproduct on $A$. A {\it counit} is a linear map $\varepsilon: A\to \Bbb C$ so that 
$$\align (\varepsilon\ot\iota)(\Delta(a)(1\ot b))&=ab \\
         (\iota\ot\varepsilon)((a\ot 1)\Delta(b))&=ab 
\endalign$$
for all $a,b\in A$.
\hfill$\square$\einspr

If $A$ has an identity and if $\Delta(A)\subseteq A\ot A$, this is equivalent with the usual conditions 
$$(\varepsilon\ot\iota)\Delta(a)=a \qquad\text{and} \qquad  (\iota\ot\varepsilon)\Delta(a)=a $$
for all $a\in A$. 
\snl
If $A$ is regular, the two conditions are equivalent with 
$$\align (\varepsilon\ot\iota)((1\ot b)\Delta(a))&=ba \\
         (\iota\ot\varepsilon)(\Delta(b)(a\ot 1))&=ba 
\endalign$$
for all $a,b\in A$. In other words, $\varepsilon$ is also a counit for the coproduct $\Delta^{\text{cop}}$.
\snl
In the $^*$-algebra case, when $\varepsilon$ is a counit, also $\overline\varepsilon$, defined by $
\overline\varepsilon(a)=\varepsilon(a^*)⁻$ will be a counit. And by taking $\frac12(\overline\varepsilon+\varepsilon)$, we see that we can assume in this case that the counit is self-adjoint, that is that $\overline\varepsilon=\varepsilon$.
\snl
There is {\it a problem} however with this definition of a counit in this context as we explain in the following remark. 

\inspr{1.9} Remark \rm 
An attempt to prove uniqueness, although trivial in the unital case, will fail. 
\snl
If the counit is an algebra map, an argument can be given as follows. Assume that $\varepsilon$ and $\varepsilon'$ are algebra maps and satisfy the axioms for a counit. Take $a,b,c\in A$ and apply $\varepsilon\ot\varepsilon'$ on the expression $(c\ot 1)\Delta(a)(1\ot b)$. Using first that $\varepsilon$ is an algebra map we get
$$\varepsilon(c)(\varepsilon'((\varepsilon\ot \iota)\Delta(a)(1\ot b)))$$
and this is, using that $\varepsilon$ is a counit,  equal to $\varepsilon(c)\varepsilon'(ab)$. Similarly, if we use first that $\varepsilon'$ is an algebra map and then that it is a counit, we find that this expression is equal to $\varepsilon'(b)\varepsilon(ca)$. This implies that 
$$\varepsilon(c)\varepsilon'(a)\varepsilon'(b)=\varepsilon'(b)\varepsilon(c)\varepsilon(a)$$
and this will imply that $\varepsilon=\varepsilon'$.
\snl
However, we know from the examples (see further) and from the theory of weak Hopf algebras, that {\it we can not expect the counit to be an algebra map}. Therefore, the above argument will  not work to prove uniqueness of the counit as defined in Definition 1.8.

\hfill$\square$\einspr

It turns out that we will need another condition on the coproduct to have that the counit, as defined in 1.8, is unique.

\iinspr{1.10} Definition \rm
A comultiplication $\Delta$ is called {\it full} if the smallest subspaces $V,W$ of $A$ satisfying 
$$\Delta(A)(1\ot A)\subseteq V \ot A \qquad \text{and} \qquad (A\ot 1)\Delta(A)\subseteq A\ot W$$
equal $A$ itself.
\hfill$\square$\einspr

We have the following easy lemma.

\iinspr{1.11} Lemma \rm
If $\Delta$ is full, the span of elements in $A$ of the form $(\omega\ot\iota)((c\ot 1)\Delta(b))$ where $b,c\in A$ and $\omega\in A'$, the linear dual space of $A$, equals $A$. Similarly for the span of elements of the form $(\iota\ot \omega)(\Delta(b)(1\ot c))$. Also conversely, if these conditions are satisfied, then the coproduct is full.

\snl\bf Proof\rm:
Suppose that such elements do not span all of $A$. Then there exists a non-zero linear functional $\varphi$ on $A$ that vanishes on all such elements. This implies that 
$$\omega((\iota\ot\varphi)((c\ot 1)\Delta(b))=0$$
for all $\omega\in A'$ and all $b,c\in A$. Then $(\iota\ot\varphi)((c\ot 1)\Delta(b))=0$ for all $b,c\in A$. If now $W$ is the kernel of $\varphi$, we see that $(c\ot 1)\Delta(b)\in A\ot W$ for all $b,c\in A$. By the assumption, this will imply that $W=A$ and therefore that $\varphi=0$. This gives a contradiction. Similarly for the other statement.
\snl
Conversely, assume that the span of these elements is all of $A$ and assume e.g.\ on the other hand that $(A\ot 1)\Delta(A)\subseteq A\ot W$ for a proper subspace $W$ of $A$. Then $(\omega(c\,\cdot\,) \ot\iota)\Delta(a)\in W$ for all $a,b,c\in A$ and all $\omega\in A'$. This will give a contradiction.

\hfill $\square$
\einspr 

One can show quite easily, using the non-degeneracy of the product, that a regular coproduct $\Delta$ is full if and only if $\Delta^{\text{cop}}$ is full.
\snl 
The following is now an easy consequence.

\iinspr{1.12} Proposition \rm 
If the coproduct is full and if there exists a counit, then it is unique.

\snl\bf Proof\rm: If $\varepsilon$ is a counit, we must have e.g.\ 
$$\varepsilon((\omega(c\,\cdot\,)\ot\iota)\Delta(b))=\omega(cb)$$ 
whenever $b,c\in A$ and $\omega\in A'$. And by the lemma, every element in $A$ has such a form. Therefore, the coproduct is uniquely determined by this formula.
\hfill$\square$\einspr

Remark that in the $^*$-algebra case, by uniqueness of the counit $\varepsilon$, we have that $\varepsilon(a^*)=\overline\varepsilon(a)$.
\snl
Of course, when the algebra has an identity, the comultiplication will be full if there is a counit. This will also be true if $A$ has no identity and if the counit is a homomorphism. However, here we do not have an identity and the counit will not be a homomorphism. Therefore, we do not automatically have that the coproduct is full.  So {\it it seems necessary} to assume both the existence of a counit as well as that the coproduct is full. 
\snl
It is also possible to formulate a weaker useful condition. One can require e.g.\ that, given $b\in A$, then $b=0$ if $\Delta(a)(1\ot b)=0$ for all $a\in A$. Similarly one can require that given $c\in A$, then $c=0$ if $(c\ot 1)\Delta(a)=0$ for all $a\in A$. This will be true if there is a counit (by the non-degeneracy of the product) or if the coproduct is full. 
\snl 
In the sequel however, for convenience, {\it we will always assume that the coproduct is full and that there is a counit}. We refer to Section 5 where we discuss possible further research on this topic.
\nl
It should be mentioned that the conditions on the coproduct we have in this section, are {\it not self-dual conditions}. In fact, it is not so obvious to see in the case of a dual pair (as defined in Definition 1.5) what these conditions on one algebra will imply on the other one. But as the conditions are reasonable enough, we will not bother about this problem at this moment. 
\snl
Also some problems remain. We will discuss more about this in Section 5 where we draw some conclusions and suggest further research.
\nl
Let us now  rather verify that these conditions are fulfilled in the examples we want to consider.

\iinspr{1.13} Example \rm
i) Let $G$ be a groupoid and consider the algebra $A=K(G)$ with the natural coproduct as given in Proposition 1.6. The counit is given by the formula $\varepsilon(f)=\sum f(e)$ where the sum is taken over all the units $G$. To show that $\Delta$ is full, take any element $p\in G$ and let $e$ be the range of $p$. Take $f=\delta_p$, the function that is $1$ in $p$ and $0$ everywhere else. Then $\Delta(f)(e,\,\cdot\,)=f$ and this shows that the right leg of $\Delta$ is all of $A$. One could also argue that $(\delta_e\ot 1)\Delta(\delta_p)=\delta_e\ot\delta_p$ for all $p$. Similarly, by taking the source of an element, we get that the left leg of $\Delta$ is all of $A$.
\snl
ii) Again let $G$ be a groupoid and consider now the algebra $B=\Bbb C G$ with the natural coproduct as given in Proposition 1.7. The counit is given by $\varepsilon(\lambda_p)=1$ for all $p\in G$. As in this case, $\Delta(B)\subseteq B\ot B$, the coproduct is automatically full.
\hfill$\square$\einspr

The situation in the case of a weak Hopf algebra is essentially trivial.

\iinspr{1.14} Example \rm
If $(A,\Delta)$ is a weak Hopf algebra, a counit exists by assumption and again the coproduct is full because already $\Delta(A)\subseteq A\ot A$.
\hfill$\square$\einspr
\nl\nl



\bf 2. The antipodes \rm
\nl
In this section, we start with a pair $(A,\Delta)$ of a non-degenerate algebra $A$ and a coproduct $\Delta$ on $A$ satisfying the assumptions discussed in Section 1. Also here, there is no need to assume that the algebra is idempotent. However, we do assume that $\Delta$ is full (cf.\ Definition 1.10 in the previous section) and that there is a unique counit $\varepsilon$ satisfying the conditions of Definition 1.8. 
\snl
We consider the maps $T_1$ and $T_2$ (as defined in Notation 1.2) and in the case of a regular coproduct also the maps $T_3$ and $T_4$ (as defined in Notation 1.4).
\snl
In the case of a multiplier Hopf algebra, the maps $T_1$ and $T_2$ are assumed to be bijective from $A\ot A$ to itself. Then the inverses $T_1^{-1}$ and $T_2^{-1}$ are given in terms of the unique antipode $S$ by the following formulas:
$$\align T_1^{-1}(a\ot b)&=((\iota\ot S)\Delta(a))(1\ot b)\\
         T_2^{-1}(a\ot b)&=(a\ot 1)((S\ot\iota)\Delta(b))
\endalign$$
for all $a,b\in A$. In fact, $S$ is constructed from the inverse maps $T_1^{-1}$ and $T_2^{-1}$ using these formulas. Recall that also the counit is constructed in the theory of multiplier Hopf algebras.
\snl
There are some remarks to be made about these formulas. We will do this later in this section (after the proof of Proposition 2.4).
\nl
\it Generalized inverses of the canonical maps \rm
\nl
Now, we no longer assume that the maps $T_1$ and $T_2$ are bijective and so we do not have the inverse maps. Nevertheless, {\it we will proceed very much as in the case of multiplier Hopf algebras} using the concept of a generalized (or von Neumann regular) inverse (see e.g.\ [G]). It can be defined in any set with an associative multiplication. In the case of our maps, we introduce it as follows. For the moment, we treat only $T_1$, but of course, similar properties can be obtained about $T_2$ and in the regular case also about $T_3$ and $T_4$.

\inspr{2.1} Definition \rm
A {\it generalized inverse} for the map $T_1$ is a linear map $R_1: A\ot A \to A \ot A$ such that 
$$  T_1 R_1 T_1 = T_1 \qquad \text{and} \qquad  R_1 T_1 R_1 = R_1.$$

\vskip -0.4cm\hfill$\square$\einspr
Multiplication is composition of maps. If we have such an inverse $R_1$ for $T_1$ and if we let $P=T_1R_1$ and $Q=R_1T_1$, then $P$ and $Q$ are idempotents. It is clear that $P$ projects onto the range of $T_1$ and that $Q$ projects onto the range of $R_1$. On the other hand, $1 - Q$ (where we use $1$ for the identity map here) is also an idempotent and it projects onto the kernel of $T_1$. Remark that, given $T_1$, the generalized inverse $R_1$ is completely determined by $P$ and $Q$. Therefore, properties of $R_1$ are expressible in terms of $P$ and $Q$. We will use this  at various places in the paper, especially in Section 3.
\nl
We will need the existence of such a generalized inverse $R_1$ for $T_1$ satisfying natural properties that are derived from related properties of $T_1$. This takes us to the following, natural pair of conditions. Recall the convention explained in the introduction. The following conditions are not assumptions in the sense that when we assume the conditions, we explicitly will mention this. 

\inspr{2.2} Condition \rm
We consider the following two conditions for the generalized inverse $R_1$ of $T_1$: 
\snl
i)\ $R_1(\iota\ot m)=(\iota\ot m)(R_1\ot \iota)$ \ on \ $A\ot A\ot A$,\newline
ii)\ $(\Delta\ot \iota)R_1=(\iota\ot R_1)(\Delta\ot\iota)$ \ on \ $A\ot A$,
\snl
where $m$ denotes multiplication, seen as a linear map from $A\ot A$ to $A$. 
\hfill$\square$\einspr

The first condition means that
$$R_1(a\ot bb')=(R_1(a\ot b))(1\ot b')$$
for all $a,b,b'\in A$. This is a natural condition as it is also satisfied by $T_1$ itself. To see where the second one comes from, multiply with an element of $A$ in the first factor from the left. Then the condition reads as
$$(T_2\ot\iota)(\iota\ot R_1)=(\iota\ot R_1)(T_2 \ot \iota)$$
and again this is natural to require because we also have coassociativity:
$$(T_2\ot\iota)(\iota\ot T_1)=(\iota\ot T_1)(T_2 \ot \iota).\tag"(2.1)"$$
Later, in Section 3, we will formulate conditions on the projection maps $T_1R_1$ and $R_1T_1$ that will imply the above formulas for $R_1$.
\snl
In a completely similar way, we can require the existence of generalized inverses for the other maps $T_2$, $T_3$ and $T_4$, satisfying similar conditions. We will not formulate them explicitly here as we continue first to focus on the maps $T_1$ and $R_1$.
\snl
At the end of this section, we will see that such  inverses  naturally exist for the examples we consider.
\snl
{\it Observe the relation} between the two conditions i) and ii) in 2.2. This is explained with duality as we see in the following remark.

\inspr{2.3} Remark \rm
Suppose that we have two algebras $A$ and $B$ with a coproduct and a non-degenerate pairing $\langle\,\cdot\, ,\,\cdot\,\rangle$ on $A\times B$ as defined in Definition 1.5. Recall the notational convention mentioned earlier. We use the same symbols for the canonical maps on $A$ and on $B$.
\snl
Consider the adjoint of the coassociativity condition (2.1). It gives 
$$(\iota\ot T_2)(T_1\ot \iota)=(T_1\ot \iota)(\iota\ot T_2). \tag"(2.2)"$$
If we then apply the counit on the third factor, we find, because $(\iota\ot\varepsilon)T_2=m$, that
$$(\iota\ot m)(T_1\ot \iota)=T_1(\iota\ot m).\tag"(2.3)"$$
Moreover, it is not hard to see that also (2.3) will imply (2.2) so that these two conditions are in fact the same.
\snl
On the other hand, if $R_1$ is a generalized inverse of $T_1$, the adjoint $R_2$ will provide a generalized inverse of $T_2$. Therefore, we have that the conditions i) and ii) in Condition 2.2 above, are in some sense dual to each other. So, from this point of view, it is also natural to assume one if the other is assumed.
\hfill$\square$\einspr
\nl
\it The antipode $S_1$\rm
\nl
We first show the existence of an antipode $S_1$, relative to $T_1$ and $R_1$. In a next item, we will consider the other cases.

\inspr{2.4} Proposition \rm
Assume that there is a generalized inverse $R_1$ for $T_1$, satisfying Condition 2.2. Then there is a linear map $S_1$ from $A$ to $L(A)$, the space of left multipliers of $A$, such that 
$$R_1(a\ot b)=\sum_{(a)} a_{(1)} \ot S_1(a_{(2)})b$$
for all $a,b\in A$.
\hfill$\square$\einspr

The last formula is given a meaning by multiplying with an element of $A$ in the first factor and from the left. The formula is completely similar to the one we encounter in Hopf algebras and multiplier Hopf algebras. We will give more comments on this result later, after the proof.

\inspr{} Proof\rm: Take $a\in A$ and define $S_1(a)$ by
$$S_1(a)b=(\varepsilon\ot\iota)R_1(a\ot b)$$
for $b\in A$. By the first condition about $R_1$, we see that indeed this formula defines $S_1(a)$ as a left multiplier of $A$ (and it justifies the notation). 
\snl
Next let $a,b,c\in A$. Then we have, using the Sweedler notation and the second condition about $R_1$ in 2.2 that
$$\align \sum_{(a)} ca_{(1)} \ot S_1(a_{(2)})b
 	&=\sum_{(a)}(\iota\ot\varepsilon\ot\iota)(\iota\ot R_1) (ca_{(1)} \ot a_{(2)} \ot b) \\
 	&=(\iota\ot\varepsilon\ot\iota)(\iota\ot R_1) (T_2\ot \iota)(c\ot a \ot b) \\
 	&=(\iota\ot\varepsilon\ot\iota) (T_2\ot \iota)(\iota\ot R_1)(c\ot a \ot b) \\
 	&=(m \ot \iota)(\iota\ot R_1)(c\ot a \ot b) \\
	&=(c\ot 1)(R_1(a \ot b)).
\endalign$$ 
And this is precisely the formula of the proposition, interpreted as above.
\hfill$\square$\einspr

We can argue that $S_1$ is completely determined by the formula above, but for that we need to use that $\Delta$ is full. We cannot conclude it from the existence of the counit.
\snl
One also has {\it some converse result}. If a map $S_1:A\to L(A)$ exist so that $R_1$ is given by the formula in the proposition, then the conditions in 2.2 are automatically satisfied. 
\snl
All this indicates that these conditions are  natural ones.
\snl
There is also {\it the following remark} to be made. Whereas by assumption, we have $R_1(a\ot b)\in A\ot A$, this is not obvious for the right hand side $\sum_{(a)} a_{(1)} \ot S_1(a_{(2)})b$ in the equation. The same phenomenon occurs in the theory of multiplier Hopf algebras. However, as we are usually working with regular ones, so that the antipode is bijective, this expression is written as $(\iota \ot S_1)((1\ot c)\Delta(a))$ where $c=S_1^{-1}(b)$ and then it is clear that this element is in $A\ot A$.
\nl
The following two well-known formulas in the theory of Hopf algebras, are also true here, but they need the correct interpretation. We will first formulate the result and prove it. Then we will indicate how the formulas are given a meaning.

\inspr{2.5} Proposition \rm
With the assumptions and the definition of $S_1$ as in the previous proposition, we have
$$ \sum_{(a)} a_{(1)}S_1(a_{(2)})a_{(3)}=a
\quad\qquad\text{and}\qquad\quad
         \sum_{(a)} S_1(a_{(1)})a_{(2)}S_1(a_{(3)})=S_1(a)$$
for all $a\in A$.

\snl\bf Proof\rm:
These two equations follow if we insert the formula for $R_1$ in terms of $S_1$, given in the previous proposition, in the formulas $T_1R_1T_1=T_1$ and $R_1T_1R_1=R_1$, and if we then apply $\varepsilon$ (or use the fullness of $\Delta$).
\hfill $\square$\einspr

\inspr{2.6} Remark \rm
i) First consider the first formula in the above proposition. Take $b\in A$. We have that $\Delta(a)(1\ot b)$ is in $A\ot A$. Write this element as $\sum_i p_i \ot q_i$ with $p_i,q_i\in A$. Then we have
$$\sum_{(a)} a_{(1)}S_1(a_{(2)})a_{(3)}b=\sum_{i,(p_i)}p_{i(1)}S_1(p_{i(2)})q_i$$
and this is well-defined in $A$ as $\sum_{(p)}p_{(1)}\ot S_1(p_{(2)})q$ is in $A\ot A$ for all $p,q\in A$.
\snl
ii) For the second formula, take again $b\in A$ and use first that $\sum_{(a)}a_{(1)}\ot S_1(a_{(2)})b$ is in $A\ot A$. If we now write this element as $\sum_i p_i \ot q_i$ with $p_i,q_i\in A$, we find
$$\sum_{(a)} S_1(a_{(1)})a_{(2)}S_1(a_{(3)})b=\sum_{i,(p_i)}S_1(p_{i(1)})p_{i(2)}q_i$$
and again this is well-defined in $A$.
\hfill$\square$\einspr

The formulas can be written as $\iota * S_1 * \iota= \iota$ and $S_1 * \iota * S_1= S_1$ in the convolution algebra of linear maps from $A$ to itself, but of course, it would be more difficult to make these formulas meaningful. Nevertheless, it is good to keep this interpretation in mind.
\snl
Needless to say that a map $S_1$ satisfying the formulas of the previous proposition, yields a generalized inverse $R_1$ of
$T_1$, satisfying the natural assumptions.

\inspr{2.7} Definition \rm We call $S_1$ {\it the antipode} relative to the inverse $R_1$ of $T_1$.
\hfill$\square$\einspr

We will also illustrate these results in the case of our examples at the end of this section.
 
\nl
\it The other antipode $S_2$ and the inverse antipodes $S_3$ and $S_4$ \rm
\nl
In a completely similar way, we can treat $T_2$. The natural assumptions for the generalized inverse $R_2$ are 
\snl
iii) $R_2(a'a\ot b)=(a'\ot 1)R_2(a\ot b)$, \newline
iv) $(\iota\ot\Delta)R_2(a\ot b)=(R_2\ot\iota)(\iota\ot\Delta)(a\ot b)$
\snl
whenever $a,a'$ and $b$ are in $A$. Again, the last equation means
$$(\iota\ot T_1)(R_2\ot\iota)=(R_2\ot\iota)(\iota\ot T_1).$$
This will yield a linear map $S_2$ from $A$ to $R(A)$, the right multipliers of $A$, satisfying and characterized by the formula
$$R_2(a\ot b)=\sum_{(b)}aS_2(b_{(1)})\ot b_{(2)}$$
for all $a,b\in A$. It is called the antipode relative to the inverse $R_2$ of $T_2$ and it also satisfies the equalities
$$\sum_{(a)} a_{(1)}S_2(a_{(2)})a_{(3)}=a 
	\quad\qquad\text{and}\qquad\quad
         \sum_{(a)} S_2(a_{(1)})a_{(2)}S_2(a_{(3)})=S_2(a)
$$
for all $a\in A$, now interpreted in $R(A)$ (that is by multiplying with an element of $A$ from the left).
\snl
In the regular case, we also have the (inverse) antipodes $S_3$ and $S_4$, relative to $T_3$ and $T_4$ respectively. They are given by the formulas
$$\align R_3(a\ot b)&=\sum_{(a)}a_{(1)}\ot bS_3(a_{(2)}) \\
         R_4(a\ot b)&=\sum_{(b)}S_4(b_{(1)})a\ot b_{(2)}
\endalign$$ 
for all $a,b\in A$. Now we have
$$\sum_{(a)} a_{(3)}S_i(a_{(2)})a_{(1)}=a 
\quad\qquad\text{and}\qquad\quad
         \sum_{(a)} S_i(a_{(3)})a_{(2)}S_i(a_{(1)})=S_i(a)
$$
for all $a\in A$ and $i=3,4$.

\inspr{2.8} Remark \rm
i) In the general case, we expect that actually $S_1(a), S_2(a)\in M(A)$ for all $a$ and that in fact $S_1=S_2$. We will then denote it by $S$ and call it the antipode. Moreover, in the regular case we expect that this antipode $S$ maps $A$ into $A$ itself, that it is bijective and that also $S_3=S_4=S^{-1}$. For this reason, we call these {\it inverse} antipodes.
\snl
ii) However, it is  clear that extra assumptions are needed. Indeed, there is no reason why these equalities should be true without any extra relations between the inverses $R_i$ themselves. On the other hand, the assumptions on these inverses seem already  quite strong so that, at least in principle, it is possible that they determine the inverses uniquely. However, this is not what we believe.
\snl
iii) In the $^*$-algebra case, it is also expected that the map $a\mapsto S(a)^*$ is involutive for the antipode and one may wonder if this would be automatic here. Again we believe that this is not the case and that extra assumptions are needed, even for this natural property.
\hfill$\square$\einspr
 
In the next item, we look at possible relations among these various antipodes and this in connection with possible properties. The results are, for the moment, {\it purely informative} in the sense that we will not rely on them further in the development. They should rather be seen as part of the motivation of what has to come. It should also help the reader to get more insight in the problems that arise. A similar {\it strategy} will be used in the next sections.
\nl
\it The antipodes, relations and properties \rm
\nl
The easiest case to consider is that of a $^*$-algebra. Indeed, when $A$ is a $^*$-algebra and $\Delta$ a $^*$-homomorphism, we have 
$$((\Delta(a)(1\ot b))^*=(1\ot b^*)\Delta(a^*)$$
for all $a,b\in A$. Therefore, from the very definitions of $T_1$ and $T_3$ we have
$$T_3=(^* \ot ^*)\circ T_1 \circ (^* \ot ^*).$$
In fact, remember that $T_3$ was defined in such a way that this equation would hold. See a remark following Notation 1.4.
\snl
Then we arrive naturally at the following result.

\inspr{2.9} Proposition \rm
Assume that $A$ is a $^*$-algebra and that $\Delta$ is a $^*$-homomorphism. Let $R_1$ be a generalized inverse for $T_1$ (satisfying the conditions in 2.2) with associated antipode $S_1$. Define $R_3$ by
$$R_3=(^* \ot ^*)\circ R_1 \circ (^* \ot ^*).$$
Then it is a generalized inverse for $T_3$ and the related antipode $S_3$ satisfies
$$S_3(a)=S_1(a^*)^*$$
for all $a\in A$.

\snl\bf Proof\rm:
The proof of the first statement is obvious. 
\snl 
To show the relation of the antipodes, take $a,b\in A$ and remark that on the one hand (because of the assumption on $R_3$)
$$ R_3(a\ot b)=\sum_{(a)}(a_{(1)}^* \ot  S_1(a_{(2)}^*)b^*)^* 
                    =\sum_{(a)}a_{(1)}\ot  bS_1(a_{(2)}^*)^*
$$
while on the other hand (because of the definition of $S_3$)
$$R_3(a\ot b)=\sum_{(a)}a_{(1)}\ot b S_3(a_{(2)}).$$
If we apply $\varepsilon$ on the first factor, we get the result.
\hfill $\square$\einspr

Of course, a similar result will be true for the triple $(T_4, R_4, S_4)$ in relation with $(T_2, R_2, S_2)$. Indeed, by convention, we also have
$$T_4=(^* \ot ^*)\circ T_2 \circ (^* \ot ^*)$$
and so, if $R_2$ is a generalized inverse of $T_2$ with associated antipode $S_2$, we get a generalized inverse $R_4$ with an associated antipode $S_4$ such that $S_4(a)=S_2(a^*)^*$ for all $a\in A$.
\snl
Compare these results with the remarks i) and iii) in 2.8.
\nl
We know from the theory of finite-dimensional weak Hopf algebras, as well as from the theory of multiplier Hopf algebras, that we will have, in the regular case, an antipode $S$ that is a bijective map from $A$ to itself and that it is both an anti-algebra and an anti-coalgebra map. 
\snl
Now, let us {\it assume that we have such a map}. So, assume for a moment that there is a bijective linear map $S:A\to A$ so that $S(ab)=S(b)S(a)$ for all $a,b$ and
$$\Delta(S(a))=\sigma (S\ot S)\Delta(a)$$
for all $a\in A$ (where we use $\sigma$ for the flip map on $A\ot A$, extended to $M(A\ot A)$).
\snl
Then we find for all $a,b$ that
$$\align T_1(S(a)\ot S(b)) &= \Delta(S(a))(1\ot S(b)) \\
                           &=\sigma(S\ot S)((b\ot 1)\Delta(a))
\endalign$$
and we see that 
$$T_2=\sigma (S^{-1}\ot S^{-1}) T_1 (S\ot S) \sigma.$$
In this case, we get the following result.

\iinspr{2.10} Proposition \rm
Assume that $S:A\to A$ is as above. Let $R_1$ be a generalized inverse of $T_1$ with associated antipode $S_1$. If we define $R_2$ by
$$R_2=\sigma (S^{-1}\ot S^{-1}) R_1 (S\ot S) \sigma,$$ 
then it is a generalized inverse for $T_2$ and the associated antipode $S_2$ satisfies 
$$S_2=S^{-1}S_1 S.$$  

\vskip -0.5 cm\hfill$\square$\einspr

The proof is again straightforward.
\snl
We have the following interesting consequence. If $S_1$ would be known to be bijective, and both an anti-algebra and an anti-coalgebra map, then we can apply the above result with $S_1$ and so we can choose $R_2$, given $R_1$, so that $S_2=S_1$. Both properties are expected.
\snl
The following two results are also related to this remark.

\iinspr{2.11} Proposition \rm
i) Let $R_1$ be a generalized inverse for $T_1$ and assume that the related antipode $S_1$ is a bijective anti-algebra map from $A$ to itself. Then $T_3$ maps $A\ot A$ to itself and if we define $R_3$ by
$$R_3=(\iota\ot S_1^{-1})T_1(\iota\ot S_1),$$
we get a generalized inverse for $T_3$. The associated antipode satisfies $S_3=S_1^{-1}$.\newline
ii) Let $R_1$ be a generalized inverse for $T_1$ and assume that the related antipode $S_1$ is a bijective anti-coalgebra map from $A$ to itself. Then $T_4$ maps $A\ot A$ to itself and if we define $R_4$ by
$$R_4=\sigma (S_1 \ot \iota)T_1(S_1^{-1}\ot \iota)\sigma,$$
we get a generalized inverse for $T_4$. For the related antipode $S_4$ we find $S_4=S_1^{-1}$.

\snl\bf Proof\rm:
i) Given $a,b\in A$, we have
$$R_1(a\ot S_1(b))=\sum_{(a)}a_{(1)}\ot S_1(a_{(2)})S_1(b)
=\sum_{(a)}a_{(1)}\ot S_1(ba_{(2)})$$
and we see that 
$$T_3=(\iota\ot S_1^{-1})R_1(\iota\ot S_1).$$
It follows that $T_3$ maps $A\ot A$ to itself. And clearly, if we put
$$R_3=(\iota\ot S_1^{-1})T_1(\iota\ot S_1),$$
we will get a generalized inverse for $T_3$. A simple calculation gives that $S_3=S_1^{-1}$ for the related antipode.
\snl
ii) Given again $a,b\in A$, we have now
$$R_1(S_1^{-1}(a)\ot b)=\sum_{(a)}S_1^{-1}(a_{(2)})\ot a_{(1)} b
=\sigma\sum_{(a)}a_{(1)}b \ot S_1^{-1}(a_{(2)})$$
and this implies that
$$T_4= \sigma(S_1\ot\iota)R_1(S_1^{-1}\ot\iota)\sigma.$$
Now, it follows that $T_4$ maps $A\ot A$ to itself and we can define a generalized inverse $R_4$ by
$$R_4=\sigma (S_1 \ot \iota)T_1(S_1^{-1}\ot \iota)\sigma$$
with associated antipode $S_4$ satisfying $S_4=S_1^{-1}$.
\hfill $\square$\einspr

If we combine the previous results we arrive at the following.

\iinspr{2.12} Proposition \rm Assume that $R_1$ is a generalized inverse for $T_1$ and that its related antipode $S_1$ is bijective from $A$ to itself and both an anti-algebra and an anti-coalgebra map. Then automatically, the coproduct is regular. And  we can define generalized inverses $R_2$, $R_3$ and $R_4$ for $T_2$, $T_3$ and $T_4$ respectively so that for the related antipodes we get
$$S_2=S_1 \qquad\qquad \text{and} \qquad\qquad S_4=S_3=S_1^{-1}.$$
\vskip -0.5 cm
\hfill$\square$\einspr

The above results are only important for motivational reasons. Indeed, we want to find conditions on the pair $(A,\Delta)$ so that we actually can prove that the assumptions about $R_1$ and $S_1$ in Proposition 2.12 hold and so we get a good antipode $S$. This is what we will do in the next section.
\nl
On the other hand, all these intermediate results motivate the following, {\it also intermediate definition}.

\iinspr{2.13} Definition \rm Let $(A,\Delta)$ be a pair of an algebra $A$ with a regular coproduct $\Delta$. Assume that the coproduct is full and that there is a counit $\varepsilon$. Let there be a bijective linear map $S:A\to A$ that is an anti-algebra and an anti-coalgebra map, satisfying 
$$\align \sum_{(a)} a_{(1)}S(a_{(2)})a_{(3)}&=a \\
         \sum_{(a)} S(a_{(1)})a_{(2)}S(a_{(3)})&=S(a)
\endalign$$
for all $a\in A$. Then we call $(A,\Delta,S)$  a {\it unifying multiplier Hopf algebra}. If moreover $A$ is a $^*$-algebra, $\Delta$ a $^*$-homomorphism and if $S$ satisfies $S(S(a)^*)^*=a$ for all $a$, we call it a unifying multiplier Hopf $^*$-algebra.
\hfill$\square$\einspr

If the algebra $A$ has an identity so that $\Delta$ maps $A$ to $A\ot A$, we speak of a unifying Hopf ($^*$-) algebra.
\snl
The conditions for the antipode are as in Proposition 2.5 and they give generalized inverses for the four canonical maps, satisfying the natural requirements (as in Condition 2.2). Also observe that in this definition, all the conditions are 'self-dual' in the sense explained in Remark 2.3.
\snl
For a unifying multiplier Hopf algebra as defined above, we can consider already the {\it source} and {\it target maps} $\varepsilon_s$ and $\varepsilon_t$, defined as
$$\varepsilon_s(a)=\sum_{(a)}S(a_{(1)})a_{(2)}
\qquad\text{and}\qquad\quad
\varepsilon_t(a)=\sum_{(a)}a_{(1)}S(a_{(2)})$$
for $a\in A$. These maps have ranges in the multiplier algebra $M(A)$. We will study these maps in greater detail in [VD-W3].
\snl
We do not plan to study this concept further. Some research on unifying (multiplier) Hopf algebras can be found in [VD-W1]. It is an earlier (unpublished) version with some of the material studied in this paper.
\nl
Of course, we can not expect that this is the final definition for weak multiplier Hopf algebras. Indeed, it is not possible to show e.g.\ that a finite-dimensional unifying (multiplier) Hopf algebra will automatically be a weak Hopf algebra. We know that we must impose extra conditions on the counit (with respect to the product of the algebra) and (dually) on the idempotent that has to replace $\Delta(1)$. This will be done in the next section. 
\nl
Before we come to that, we look at examples to illustrate the definitions and results in this section.
\nl
\it Examples \rm
\nl
We briefly consider the two examples coming from a groupoid, as well as the case of an ordinary weak Hopf algebra.
\iinspr{2.14} Example \rm
i) Consider first the example of Proposition 1.6. There exists an antipode $S$, defined in the usual way by $S(f)(p)=f(p^{-1})$ for $f\in K(G)$ and $p\in G$. The basic property to use is that $pp^{-1}p=p$ as well as $p^{-1}pp^{-1}=p^{-1}$ for all elements $p$ in $G$. We obviously get all the conditions to have a unifying multiplier Hopf algebra as in Definition 2.13.
\snl
ii) Next, consider the example of Proposition 1.7. Here an antipode $S$ is given by the formula $S(\lambda_p)=\lambda_{p^{-1}}$ for all $p$ in $G$. We are then exactly in the same situation as in the previous example.
\hfill$\square$\einspr

\iinspr{2.15} Example \rm
Also in a weak Hopf algebra, all the conditions of Definition 2.13 are fulfilled (see  Lemma  2.8 and Proposition 2.10 in [B-N-S] and  Proposition 2.3.1 in [N-V2]). So any weak Hopf algebra will be a unifying Hopf algebra in the sense of the above definition.
\hfill$\square$\einspr

\nl



\bf 3. The idempotent $E$ in $M(A\ot A)$ and related elements \rm
\nl
In this section, we will discuss the {\it final assumptions} that eventually will yield the definition of a weak multiplier Hopf algebra in the next section. 
\snl
The starting point is again a pair $(A,\Delta)$ of a non-degenerate  algebra $A$ and a coproduct $\Delta$ on $A$ (in the sense of Definition 1.1). Later in this section, we will need that $A$ is idempotent and so we assume that this is the case. We also assume that $\Delta$ is full (see Definition 1.10) and that there is a (unique) counit $\varepsilon$ on $A$ (as in Definition 1.8). In the case of a $^*$-algebra, we assume that $\Delta$ is a $^*$-homomorphism.
\snl
We consider the maps $T_1$ and $T_2$ and now we assume that we have generalized inverses $R_1$ and $R_2$ satisfying the necessary conditions (cf.\ Condition 2.2) so that there are associated antipodes $S_1$ and $S_2$ as discussed in the previous section. If the coproduct $\Delta$ is regular (cf.\ Definition 1.3), we also consider the maps $T_3$ and $T_4$ and we assume that they have generalized inverses $R_3$ and $R_4$ with associated (inverse) antipodes $S_3$ and $S_4$. 
\snl
{\it For the moment}, we will not require any of the relations among these inverses as discussed at the end of the previous section (as in Proposition 2.9 and further). Some of them however will be considered again at the end of this section.
\snl 
{\it On the other hand}, in this section, we will be considering {\it extra conditions} on the idempotent maps $T_1R_1$ and $T_2R_2$ and on the idempotent maps $R_1T_1$ and $R_2T_2$. We will then see what the conditions on the inverses $R_1$ and $R_2$ imply in the new situation. And conversely, we will try to obtain the conditions on these inverses in turn from the conditions and formulas involving the idempotents.
\snl
We will discuss these new assumptions and give different arguments to see that they are natural. We will show that they are satisfied in our basic examples. And what is very important, we will see how they all are intimately related and that they will {\it complete the puzzle}, finally providing a good notion of a weak multiplier Hopf algebra in the next section.
\nl
\it The idempotent $E$ \rm
\nl
The first assumption is about the behavior of the idempotents $T_1R_1$ and $T_2R_2$ with respect to multiplication. It is a preliminary assumption and we will later strengthen it (see Assumption 3.4 below). We formulate it as follows.

\inspr{3.1} Assumption \rm
We assume that the generalized inverses $R_1$ and $R_2$ can be chosen in such a way that also 
$$\align T_1R_1(aa'\ot b) & = (T_1R_1(a\ot b))(a'\ot 1)\tag"(i)"\\
			T_2R_2(a\ot b'b) & = (1\ot b')T_2R_2(a\ot b)\tag"(ii)"
\endalign$$
for all $a,a',b,b'\in A$.
\hfill$\square$\einspr

Because we already have a similar property with multiplication for the other factor, both for $T_1$ and $R_1$ (Condition 2.2), if we also assume this property, it follows  that $T_1R_1$ is a left multiplier of $A\ot A$. Similarly $T_2R_2$ is a right multiplier. This makes the following notation possible. 

\inspr{3.2} Notation \rm
Define a left multiplier $E$ of $A\ot A$ and a right multiplier $E'$ of $A\ot A$ by
$$E(a\ot b)=T_1R_1(a\ot b) \qquad\qquad\text{and}\qquad\qquad (a\ot b)E'=T_2R_2(a\ot b)$$
for $a,b\in A$.
\hfill$\square$\einspr

We have $E^2=E$ and ${E'}^2=E'$ because $T_1R_1$ and $T_2R_2$ are idempotent maps. 
\snl
Also observe that $E(A\ot A)$ is the range of $T_1$ and that $(A\ot A)E'$ is the range of $T_2$. Remark however that the assumptions require not only that the ranges of $T_1$ and $T_2$ have this form but also that one can choose the inverses $R_1$ and $R_2$  with these extra properties. Of course, this choice is only possible under these assumptions on the ranges. 
\nl
The following lemma will yield the uniqueness of such idempotents under mild and natural conditions.

\inspr{3.3} Lemma \rm
Suppose that we have two idempotents $E$ and $F$ in $M(A\ot A)$ so that 
$$F(A\ot A)\subseteq E(A\ot A) \qquad\qquad\text{and}\qquad\qquad (A\ot A)E\subseteq (A\ot A)F,$$
then $E=F$.

\snl \bf Proof\rm: 
Take $a,b\in A$. Use again $1$ for the identity map in this context. From the  assumptions we get
$$\align (1-E)F(a\ot b)&\in (1-E)E(A\ot A)\\
         (a\ot b)E(1-F)&\in (A\ot A)F(1-F)
\endalign$$
and because $(1-E)E=0$ and $F(1-F)=0$, we see that $(1-E)F=0$ as well as $E(1-F)=0$. This gives $E=F$.
\hfill$\square$
\einspr

So, if we have one idempotent $E\in M(A\ot A)$ so that $E(A\ot A)=\text{Ran}(T_1)$ and \newline 
$(A\ot A)E=\text{Ran}(T_2)$, this idempotent is unique. Let us therefore replace Assumption 3.1 by the following stronger assumption.

\inspr{3.4} Assumption \rm
Assume that we have generalized inverses $R_1$ and $R_2$ such that there is an idempotent $E\in M(A\ot A)$ satisfying 
$$T_1R_1(a\ot b)=E(a\ot b) \qquad\qquad\text{and}\qquad\qquad T_2R_2(a\ot b)=(a\ot b)E$$
for all $a,b\in A$.
\hfill$\square$\einspr

Of course, this idempotent will play the role of $\Delta(1)$ in our theory. In fact, we have the following result endorsing this statement.

\inspr{3.5} Proposition \rm
The idempotent $E$ is the smallest one in $M(A\ot A)$ satisfying 
$$E\Delta(a)=\Delta(a) \qquad\qquad\text{and}\qquad\qquad \Delta(a)E=\Delta(a)$$
for all $a\in A$.

\bf \snl Proof\rm:
By the definition of $E$ we have
$$E(T_1(a\ot b))=T_1R_1T_1(a\ot b)=T_1(a\ot b)$$
and so $E(\Delta(a)(1\ot b)=\Delta(a)(1\ot b)$ for all $a,b\in A$. Because the product is non-degenerate, we get $E\Delta(a)=\Delta(a)$ for all $a$. Similarly, by using that $(a\ot b)E=T_2R_2(a\ot b)$, we get from the definition of $T_2$ that $\Delta(a)E=\Delta(a)$ for all $a$.
\snl
Now suppose that $E'$ is another idempotent in $M(A\ot A)$ so that $E'\Delta(a)=\Delta(a)$ and $\Delta(a)E'=\Delta(a)$ for all $a$. Then it will follow that $E'T_1(a\ot b)=T_1(a\ot b)$ and this will give $E'E=E$. Similarly we get $EE'=E$. This is precisely what we mean by saying that $E$ is smaller than $E'$.
\hfill$\square$
\einspr

It follows  that $E^*=E$ in the involutive case. Indeed, because now $\Delta$ is assumed to be a $^*$-homomorphism, it will follow that $E^*\Delta(a)=\Delta(a)=\Delta(a)E^*$ for all $a$ and by the previous result we have $EE^*=E^*E=E$. Taking adjoints we get $E=E^*$. 
\snl 
In the regular case, we should impose similar assumptions on the inverses $R_3$ and $R_4$ giving the same idempotent. That is, we should require then (because of our conventions) that
$$(a\ot b)E=T_3R_3(a\ot b) \qquad\qquad\text{and}\qquad\qquad E(a\ot b)=T_4R_4(a\ot b)$$
for all $a,b\in A$. We will come back to this later (see a remark after Proposition 3.9 and also Proposition 4.9 in the next section).
\nl
Before we start investigating results on this idempotent $E$, we first see what happens in the examples and we show that the extra conditions are fulfilled in the standard examples.

\inspr{3.6} Examples \rm
i) First we consider the case of Proposition 1.6, see also Example 2.14.i. If we again use the obvious antipode $S$, both for $R_1$ and $R_2$ we find that 
$$\align (T_1R_1f)(p,q)&=f(pqq^{-1},q) \\
         (T_2R_2f)(p,q)&=f(p,p^{-1}pq)
\endalign$$
when $f\in K(G)$ and when $p,q\in G$ are elements so that $pq$ is defined (that is when $t(q)=s(p)$). In the other case, we get $0$. So we see that in both cases $T_1R_1f=Ef$ and $T_2R_2f=Ef$ where $E$ is the function on $G\times G$ that is one on pairs $(p,q)$ for which $pq$ is defined and $0$ on other pairs. Remark that the algebra is abelian. So, the Assumption 3.4 is fulfilled. Of course, $E$ is the obvious candidate for $\Delta(1)$ as we observed already in a remark following Proposition 1.6.
\snl
ii) Next, consider the dual case as in Proposition 1.7, see also Example 2.14.ii. Again we consider the obvious antipode $S$ for defining $R_1$ and $R_2$. We find, for $p,q\in G$ that
$$\align T_1R_1(\lambda_p\ot\lambda_q) &= \lambda_p \ot \lambda_{pp^{-1}q} \\
         T_2R_2(\lambda_p\ot\lambda_q) &= \lambda_{pq^{-1}q} \ot \lambda_q 
\endalign$$
where in the first equation we need that $t(q)=t(p)$ and in the second equation that $s(q)=s(p)$ in order to get a non-zero outcome. So if $E$ is defined as $\sum_e \lambda_e\ot \lambda_e$ (where the sum is taken over all units $e$), it follows that $T_1R_1$ is left multiplication with $E$ and $T_2R_2$ is right multiplication with $E$. Again Assumption 3.4 is fulfilled. Also here $E$ is the obvious candidate for $\Delta(1)$ as we also observed in a remark following Proposition 1.7.
\hfill$\square$
\einspr

\inspr{3.7} Example \rm
In the case of a finite-dimensional weak Hopf algebra $(A,\Delta)$ we find the same results. If $S$ is the antipode and $R_1$ and $R_2$ are defined using this antipode, we get (using Definition 2.1.1  and Proposition 2.2.1 of [N-V2])
$$\align T_1R_1(a\ot b)&=\sum_{(a)}a_{(1)}\ot a_{(2)}S(a_{(3)})b \\
		&=\sum_{(a)}a_{(1)}\ot \varepsilon_t(a_{(2)})b\\
		&=\Delta(1)(a\ot b)
\endalign$$
for all $a,b\in A$. Similarly $T_2R_2(a\ot b)=(a\ot b)\Delta(1)$. So the assumptions are fulfilled with $E=\Delta(1)$ as expected.
\hfill$\square$
\einspr

\newpage
\it Further properties of the idempotent $E$ \rm
\nl
We have seen that the idempotent $E$ has to be considered as a replacement for $\Delta(1)$ and we have made this precise in Proposition 3.5 above.  This has a few immediate and important consequences. 

\inspr{3.8} Proposition \rm
There is a unique homomorphism $\Delta_1:M(A)\to M(A\ot A)$ that extends $\Delta$ and satisfies $\Delta_1(1)=E$.
\hfill$\square$\einspr

We refer to  Appendix A (Proposition A.2) for a detailed proof of this (and related) result(s). Just remark that this extension is characterized by
$$\align \Delta_1(m)(\Delta(a)(1\ot b)&=\Delta(ma)(1\ot b) \\
		  ((a\ot 1)\Delta(b))\Delta_1(m)&=(a\ot 1)\Delta(bm),
\endalign$$
for all $a,b\in A$  and $m\in M(A)$, and by the requirement that 
$$\Delta_1(m)=\Delta_1(m)E\qquad\qquad\text{and}\qquad\qquad \Delta_1(m)=E\Delta_1(m)$$
for all $m\in M(A)$.
\snl
It is here that we seem to need that the algebra $A$ is idempotent (see the appendix).
\snl
As usual, we will denote this extension still by $\Delta$.
\snl
Similarly, we can extend the homomorphisms $\Delta\ot\iota$ and $\iota\ot\Delta$ to homomorphisms from $M(A\ot A)$ to $M(A\ot A\ot A)$ using the requirement
$$\align (\Delta\ot\iota)(m)&=(\Delta\ot\iota)(m)(E\ot 1) \\
         (\Delta\ot\iota)(m)&=(E\ot 1)(\Delta\ot\iota)(m)        
\endalign$$
for all $m\in M(A\ot A)$ and similarly for $\iota\ot\Delta$. Again see Appendix A, Proposition A.5.
\nl
With these definitions, we have the following result.

\inspr{3.9} Proposition \rm
We have 
$$ (\Delta\ot\iota)(E)=(1\ot E)(E\ot 1) 
		\qquad\quad\text{and}\qquad\quad
         (\iota\ot\Delta)(E)=(1\ot E)(E\ot 1).
$$

\snl\bf Proof\rm:
We have that 
$$(\Delta\ot\iota)T_1R_1=(\iota\ot T_1R_1)(\Delta\ot \iota)$$
as we have these commutation rules for $T_1$ (by definition) and for $R_1$ (by assumption). If we apply this equation to $a\ot b$ and use that $\Delta\ot\iota$ is still a homomorphism on $M(A\ot A)$, we find
$$(\Delta\ot\iota)(E)(\Delta(a)\ot b)=(1\ot E)(\Delta(a)\ot b)$$
for all $a,b\in A$. This implies that
$$(\Delta\ot\iota)(E)=(\Delta\ot\iota)(E)(E\ot 1)=(1\ot E)(E\ot 1)$$
where we have used that $(\Delta\ot\iota)(m)=(\Delta\ot\iota)(m)(E\ot 1)$ for all $m\in M(A\ot A)$. This proves the first formula.
\snl
Similarly, from 
$$(\iota\ot\Delta)T_2R_2=(T_2R_2\ot\iota)(\iota\ot\Delta),$$
we get the other formula.
\hfill$\square$
\einspr

So, we see that in particular $(\Delta\ot\iota)(E)=(\iota\ot\Delta)(E)$, something that we showed already in greater generality in Proposition A.8 of the appendix. We also know (see Proposition A.7) that it follows that 
$$(\Delta\ot\iota)(\Delta(m))=(\iota\ot\Delta)(\Delta(m))$$
for every multiplier $m\in M(A)$, as expected.
\nl
We will now give some attention to the {\it regular case}, mainly for motivational reasons. 
\snl
We begin with the involutive case. So, let $A$ be a $^*$algebra and $\Delta$ a $^*$-homomorphism.  It has been shown to be  a consequence of Proposition 3.5 that $E^*=E$ in this case. Then, as we have the formula $(\iota\ot\Delta)(E)=(1\ot E)(E\ot 1)$, it follows by taking the adjoint that also $(\iota\ot\Delta)(E)=(E\ot 1)(1\ot E)$ and $1\ot E$ and $E\ot 1$ commute.
\snl
As we have seen in Proposition 2.9, given $R_1$, we can choose $R_3=(^* \ot ^*)\circ R_1 \circ (^* \ot ^*)$ and as we have a similar formula relating $T_3$ with $T_1$, we get also
$$T_3R_3=(^* \ot ^*)\circ T_1R_1 \circ (^* \ot ^*).$$
Because $E^*=E$, it follows from this that $T_3R_3(a\ot b)=(a\ot b)E$.
\snl 
This suggests that also in the general regular case, it is natural to assume that 
$T_3R_3(a\ot b)=(a\ot b)E$ for all $a,b$ and if we use the given commutation rules of $T_3$ and $R_3$ with $\Delta\ot\iota$, we get 
$$(\Delta(a)\ot b)(\Delta\ot\iota)(E)=(\Delta(a)\ot b)(\iota\ot E)$$
and this will imply that 
$$(\Delta\ot\iota)(E)=(E\ot 1)(1\ot E).$$
As we have shown already in Proposition 3.9 that also $(\Delta\ot\iota)(E)=(1\ot E)(E\ot 1)$, again we  arrive at the fact that $1\ot E$ and $E\ot 1$ commute.
\snl
This suggests the following natural assumption.

\iinspr{3.10} Assumption \rm
We assume that 
$$(E\ot 1)(1\ot E)=(1\ot E)(E\ot 1).$$
\vskip -0.5 cm
\hfill$\square$\einspr

This implies that 
$$(\Delta\ot\iota)(E)=(1\ot E)(E\ot 1)=(E\ot 1)(1\ot E).$$
\snl
From this condition, we can prove the following new commutation rules. Observe the difference with the rule in Condition 2.2.

\iinspr{3.11} Proposition \rm 
\snl
i) $(\iota\ot\Delta)T_1R_1=(T_1R_1\ot \iota)(\iota\ot\Delta)$,\newline
ii) $(\Delta\ot\iota)T_2R_2=(\iota\ot T_2R_2)(\Delta\ot\iota)$.

\snl\bf Proof\rm:
Given the fact that $T_1R_1(a\ot b)=E(a\ot b)$ for all $a,b\in A$, we see that the first formula can be written as
$$(\iota\ot\Delta)(E)(a\ot\Delta(b))=(E\ot 1)(a\ot \Delta(b))$$
for all $a,b\in A$. We can safely cancel $a$ and as $\Delta(b)=E\Delta(b)$, we see that 
i) follows from the Assumption 3.10. Similarly ii) will follow.
\hfill$\square$
\einspr

It is also possible to show some converse result. Indeed, if one of the properties in Proposition 3.11 is satisfied, then also Assumption 3.10 will follow. To see this, 
multiply the formula in the proof with $1\ot 1\ot c$ and use that the range of $T_1$ is precisely $E(A\ot A)$. This shows that i)  implies $(\Delta\ot\iota)(E)=(E\ot 1)(1\ot E)$. 
\snl
Before we verify the new conditions in the case of the examples, we make one more remark.
\snl
If we rewrite the assumption in 3.1 about $T_1R_1$ and $T_2R_2$ as
$$\align 
T_1R_1(m^{\text{op}}\ot \iota)&=(m^{\text{op}}\ot \iota)(\iota\ot T_1R_1)\\
T_2R_2(\iota\ot m^{\text{op}})&=(\iota\ot m^{\text{op}})(T_2R_2\ot\iota),\\
\endalign$$
we see a similarity with the formulas in the proposition above. The conditions are (almost) dual to each other. We would get dual conditions if they were written in terms of $R_1T_1$ and $R_2T_2$ and not for $T_1R_1$ and $T_2R_2$. We consider this in the next item where we also collect all these rules.
\snl
Now, observe that Assumption 3.10 is {\it fulfilled in the case of our examples}. In the case $A=K(G)$ (Example 3.6.i) there is no problem as the algebra is abelian. In the case $A=\Bbb C G$ (Example 3.6.ii), it is immediately seen that
$$(\iota\ot\Delta)(E)=(\Delta\ot \iota)(E)=\sum_e \lambda_e \ot \lambda_e \ot \lambda_e$$ 
(where the sum is taken over the units $e$ of $G$) and this is indeed equal to
$(1\ot E)(E\ot 1)$ and $(E\ot 1)(1\ot E)$. 
\snl
In the case of a weak Hopf algebra, these equalities are all true as they are part of the axioms. 
\nl
\it Conditions on the idempotent maps $R_1T_1$ and $R_2T_2$ \rm
\nl
In the previous item, we looked at the idempotent maps $T_1R_1$ and $T_2R_2$ and found reasonable conditions giving the idempotent $E$ in $M(A\ot A)$, playing the role of $\Delta(1)$ and satisfying some expected properties. Now, we will look at the idempotents $R_1T_1$ and $R_2T_2$ and find natural conditions by duality. This is based on the idea that the maps $R_1T_1$ and $R_2T_2$, in the case of a dual pair, are adjoint to the maps $T_2R_2$ and $T_1R_1$ respectively.
\snl 
We have the (original) properties
$$\align T_1(\iota\ot m) & = (\iota\ot m)(T_1\ot\iota) \\
         T_2(m\ot\iota)  & = (m\ot\iota)(\iota\ot T_2)
\endalign$$
on $A\ot A\ot A$ and
$$\align (\Delta\ot\iota)T_1 & = (\iota\ot T_1)(\Delta\ot\iota) \\
         (\iota\ot\Delta)T_2 & = (T_2\ot\iota)(\iota\ot\Delta)
\endalign$$
on $A\ot A$. The second set is dual to the first one. As we have assumed the same commutation rules for the generalized inverses $R_1$ and $R_2$ in Section 2, we also get these  rules for the composed maps $T_1R_1$, $R_1T_1$, $T_2R_2$ and $R_2T_2$. In what follows, we will call these the {\it original commutation rules}. 
\snl
Now, in the first two items of this section, we have added what we will call the {\it new commutation rules}. These can be formulated as
$$\align T_1R_1(m^{\text{op}}\ot\iota) & = (m^{\text{op}}\ot\iota)(\iota\ot T_1R_1) \\
         T_2R_2(\iota\ot m^{\text{op}}) & = (\iota\ot m^{\text{op}})(T_2R_2\ot\iota)
\endalign$$
on $A\ot A\ot A$ and
$$\align (\iota\ot\Delta)T_1R_1 & = (T_1R_1 \ot \iota)(\iota\ot \Delta) \\
         (\Delta\ot\iota)T_2R_2 & = (\iota\ot T_2R_2)(\Delta\ot\iota)
\endalign$$
on $A\ot A$. 
\snl
The first equations are a reformulation of the assumptions in 3.1 as we mentioned already whereas the second pair of equations follows from the assumption in 3.10 as proven in Proposition 3.11.
\snl
These new conditions however are not self-dual because the dual conditions should give commutation rules involving the other idempotent maps $R_1T_1$ and $R_2T_2$.
\nl
A straightforward application of the rules on duality (as in Definition 1.5) gives that these four new rules dualize as follows. We formulate them as a {\it second set of new commutation rules}.

\iinspr{3.12} Assumption \rm
We assume that 
$$\align R_1T_1(m\ot\iota) & = (m\ot\iota)(\iota\ot R_1T_1) \\ 
         R_2T_2(\iota\ot m) & = (\iota\ot m)(R_2T_2\ot\iota)
\endalign$$
on $A\ot A\ot A$ and also
$$\align (\iota\ot\Delta^{\text{cop}})R_1T_1 & = (R_1T_1 \ot \iota)(\iota\ot \Delta^{\text{cop}}) \\
         (\Delta^{\text{cop}}\ot\iota)R_2T_2 & = (\iota\ot R_2T_2)(\Delta^{\text{cop}}\ot\iota)
\endalign$$
on $A\ot A$.
\hfill$\square$\einspr

So, these assumptions are {\it justified using duality}. 
\snl
They are also natural from another point of view as we argue in the following remark.

\iinspr{3.13} Remark \rm
i) Suppose that the antipode $S_1$, defined relative to $R_1$ in Section 2, is bijective from $A$ to itself and an anti-homomorphism (and that $\Delta$ is regular). Then, as we have already seen before, we have 
$$R_1(a\ot S_1(b))=\sum_{(a)}a_{(1)}\ot S_1(a_{(2)})S_1(b) = \sum_{(a)}a_{(1)}\ot S_1(ba_{(2)})$$
for all $a,b\in A$ and so
$$R_1=(\iota\ot S_1)T_3(\iota\ot S_1^{-1}).$$
If (as in Proposition 2.11) we choose $R_3$ so that also
$$T_1=(\iota\ot S_1)R_3(\iota\ot S_1^{-1}),$$
we find
$$R_1T_1=(\iota\ot S_1)T_3R_3(\iota\ot S_1^{-1}).$$
Now, as explained earlier, we expect that 
$$T_3R_3(a\ot b)=(a\ot b)E$$ 
and so
$$R_1T_1(a\ot S_1(b))=(\iota\ot S_1)((a\ot b)E).$$
ii) This will justify the assumption that 
$$R_1T_1(a'a\ot b)=(a'\ot 1)R_1T_1(a\ot b)$$
for all $a,a',b\in A$.
Similarly, if we work with the antipode $S_2$ relative to $R_2$, we would obtain
$$R_2T_2(a\ot bb')=(R_2T_2(a\ot b))(1\ot b')$$
for all $a,b,b'\in A$.
\hfill$\square$\einspr

Of course the argument that uses duality is a better argument as it does not assume anything about the antipodes $S_1$ and $S_2$.
\snl
Remark again that these assumptions are in the first place conditions about the kernels $\text{Ker}(T_1)$ and $\text{Ker}(T_2)$, just as the assumptions in 3.1 are first of all conditions on the ranges $\text{Ran}(T_1)$ and $\text{Ran}(T_2)$ as mentioned already.
\nl
Observe that from the discussion in the remark, we see that we expect
$$\align 
	R_1T_1(a\ot b)&=(a\ot 1)F_1(1\ot b) \\
	R_2T_2(a\ot b) & = (a\ot 1)F_2(1\ot b)
\endalign$$
for all $a,b\in A$ where $F_1=(\iota\ot S_1)E$ and $F_2=(S_2\ot \iota)E$. In fact, as a consequence of the extra assumptions in 3.12, we can obtain such formulas, giving the idempotents $R_1T_1$ and $R_2T_2$ in terms of a multiplier. Compare with the formulas in Notation 3.2.

\iinspr{3.14} Proposition \rm
There is a right multiplier $F_1$ of $A\ot A^{\text{op}}$ and a left multiplier $F_2$ of $A^{\text{op}}\ot A$ such that
$$\align R_1T_1(a\ot b) & = (a\ot 1)F_1(1\ot b) \\
         R_2T_2(a\ot b) & = (a\ot 1)F_2(1\ot b)
\endalign$$
in $A\ot A$ for all $a,b\in A$.
\hfill$\square$\einspr

The proof is very simple. We use the original and the new module properties:
$$\align R_1T_1(a'a\ot b) & = (a'\ot 1)(R_1T_1(a\ot b)) \\
         R_1T_1(a\ot bb') & = (R_1T_1(a\ot b))(1\ot b')
\endalign$$
for all $a,a',b,b'\in A$ and similarly for $R_2T_2$.
\snl
Remark that, for the idempotents $T_1R_1$ and $T_2R_2$, we have added an extra assumption (Assumption 3.4) relating the two. We should also expect some relation here between $F_1$ and $F_2$ but this is not so simple. We will come back to this in the next item of this section.
\snl
On the other hand, from  the original and the new commutation rules with the coproduct, we get the following formulas.  Compare them with the formulas involving $E$ in Proposition 3.9. We use the leg-numbering notation as explained in the introduction. 

\iinspr{3.15} Proposition \rm
On the one hand we get
$$\align (\Delta\ot\iota)F_1 & = (E\ot 1)(1\ot F_1)\\
         (\iota\ot\Delta)F_2 & = (F_2\ot 1)(1\ot E)
\endalign$$ 
while on the other hand we get
$$\align (\iota\ot\Delta)F_1 & = (F_1)_{13}(1\ot E)\\
         (\Delta\ot\iota)F_2 & = (E\ot 1)(F_2)_{13}.
\endalign$$
\einspr

We need to be a little careful with the interpretation of these formulas. We have extended the homomorphisms $\Delta\ot\iota$ and $\iota\ot\Delta$ to the multiplier algebra $M(A)$ in a certain way. Here we are applying these maps to a left and a right multiplier. However, we can interpret these formulas by multiplying at the right place with elements of $A$ as we do in the proof below.

\iinspr{} Proof\rm: If we combine the original formulas
$$\align (\Delta\ot\iota)R_1T_1 & = (\iota\ot R_1T_1)(\Delta\ot\iota) \\
         (\iota\ot\Delta)R_2T_2 & = (\iota\ot T_2R_2)(\iota\ot\Delta)
\endalign$$
with the formulas in Proposition 3.14, we get easily the first results. Consider e.g.\ the first one. Apply it to $a\ot b$ with $a,b\in A$. We get 
$$(\Delta\ot\iota)((a\ot 1)F_1(1\ot b))=(\Delta(a)\ot\iota)(1\ot F_1)(1\ot 1\ot b).$$
If we cancel $b$, we find 
$$(\Delta\ot\iota)((a\ot 1)F_1)=\Delta(a)\ot 1)(1\ot F_1)$$
and this precisely means $(\Delta\ot\iota)F_1  = (E\ot 1)(1\ot F_1)$ by the considerations of the appendix about extending the map $\Delta\ot\iota$. Similarly for the other case.
\snl
If we combine the new formulas
$$\align (\iota\ot\Delta^{\text{cop}})R_1T_1 & = (R_1T_1 \ot \iota)(\iota\ot \Delta^{\text{cop}}) \\
         (\Delta^{\text{cop}}\ot\iota)R_2T_2 & = (\iota\ot R_2T_2)(\Delta^{\text{cop}}\ot\iota)
\endalign$$
with the formulas in Proposition 3.14, we get, after applying the flip map,
$$\align (\iota\ot\Delta)(F_1(1\ot b))&=(F_1)_{13}(1\ot \Delta(b)) \\
        (\Delta\ot\iota)((a\ot 1)F_2)&=(\Delta(a)\ot 1)(F_2)_{13}
\endalign$$
for all $a$ and $b$ in $A$. This gives, with the right interpretation, the second set of formulas in the proposition.
\hfill$\square$
\einspr 

The formulas in Proposition 3.15 are not very surprising. Indeed, if as we expect, 
$$F_1=(\iota\ot S)E
\qquad\quad\text{and}\quad\qquad
F_2=(S\ot\iota)E,$$
they essentially are the same as the formulas that give $(\Delta\ot\iota)E$ as in Proposition 3.9 and Assumption 3.10.
\snl
Before we continue our investigations, we want to make another very important remark.

\iinspr{3.16} Remark \rm
i) With the introduction of $E$, we make a choice for the range projections of $T_1$ and $T_2$. With $F_1$ and $F_2$ we do the same for the projections on the kernels of $T_1$ and $T_2$. Then the generalized inverses $R_1$ and $R_2$ are determined. So the elements $E$, $F_1$ and $F_2$ determine the inverses $R_1$ and $R_2$. See a remark following Definition 2.1.
\snl
ii) The commutation rule of the generalized inverse $R_1$ with $\iota\ot m$  (as in Condition 2.2) is  now automatically satisfied because this commutation rule is not only true for $T_1$ , but also for $T_1R_1$ and $R_1T_1$. See again the remark following Definition 2.1. Similarly for the inverse $R_2$.
\snl
iii) Finally, the commutation rule of $R_1$ with $\Delta\ot\iota$ (as in Condition 2.2) is now also automatically satisfied because this is true for $T_1$ as well as for $T_1R_1$ and $R_1T_1$ by the formulas involving $\Delta$ and $E$ (Proposition 3.9 and 3.10). Similarly for the inverse $R_2$.
\snl
iv) As a consequence, the various conditions on these idempotent maps will yield the existence of the associated antipodes $S_1$ and $S_2$ by applying the results of Section 2. Therefore, in some sense, the conditions in this section {\it override} the ones in the previous section.
\hfill$\square$\einspr

\it The relations between $E$, $F_1$ and $F_2$ \rm
\nl
For the ranges of $T_1$ and $T_2$, it was easy to relate them. We simply took $E=E'$ (from Notation 3.2) in Assumption 3.4. However, as mentioned already, there is no equivalent simple condition for the multipliers $F_1$ and $F_2$. 
\snl
We have seen from an earlier discussion in this section that if $S_1$ is a bijective anti-algebra map, we expect $F_1=(\iota \ot S_1)E$. Similarly, if $S_2$ is a bijective anti-algebra map, we expect $F_2=(S_2\ot\iota)E$. 
\snl
There is however another way to relate the multipliers $E$, $F_1$ and $F_2$ provided we assume that the antipodes $S_1$ and $S_2$ coincide. We will obtain this in the next proposition.
The result is {\it rather remarkable} and {\it very important} for the {\it further approach} and we will discuss it later. In particular, we will see how this also implies in turn that the antipode thus obtained satisfies the expected properties. 

\iinspr{3.17} Proposition \rm
The antipodes $S_1$ and $S_2$ coincide if and only if
$$E_{13}(F_1\ot 1)=E_{13}(1\ot E) \qquad\text{and}\qquad (1\ot F_2)E_{13}  =(E\ot 1)E_{13}.$$
\einspr
\vskip -0.2 cm
Recall that $S_1$ is defined as a left multiplier and that $S_2$ is a right multiplier. That $S_1$ and $S_2$ coincide simply means that $c(S_1(a)b)=(cS_2(a))b$ for all $a,b,c\in A$.

\iinspr{} Proof\rm:
Because $E(p\ot q)=T_1R_1(p\ot q)$ for all $p,q\in A$ we have 
$$\sum_{(a)}(E\ot 1)(a_{(1)}\ot 1\ot a_{(2)})=\sum_{(a)}a_{(1)}\ot a_{(2)}S_1(a_{(3)})\ot a_{(4)}$$
and because $(p\ot 1)F_2(1\ot q)=R_2T_2(p\ot q)$ for all $p,q\in A$, we get
$$\sum_{(a)}(a_{(1)}\ot 1\ot 1)(1\ot F_2)(1\ot 1\ot a_{(2)})=\sum_{(a)}a_{(1)}\ot a_{(2)}S_2(a_{(3)})\ot a_{(4)}$$
for all $a$ in $A$. Therefore we see that 
$(E\ot 1)E_{13}=(1\ot F_2)E_{13}$ if and only if 
$$\sum_{(a)}a_{(1)}S_1(a_{(2)})=\sum_{(a)}a_{(1)}S_2(a_{(2)})\tag"(3.1)"$$
for all $a$. Similarly we have that $E_{13}(1\ot E)=E_{13}(F_1\ot 1)$ if and only if
$$\sum_{(a)}S_1(a_{(1)})a_{(2)}=\sum_{(a)}S_2(a_{(1)})a_{(2)}\tag"(3.2)"$$
for all $a$. 
Remark that the left hand side of (3.1) and of (3.2) are defined as left multipliers of $A$, whereas the right hand side of (3.1) and of (3.2) are defined as right multipliers. Also here, equality is in the sense as explained above.
\snl
In particular, if $S_1=S_2$  we get the desired equations.
\snl
On the other hand, if these formulas are satisfied, we find, by the result above,
that
$$\align S_1(a)&=\sum_{(a)}S_1(a_{(1)})a_{(2)}S_1(a_{(3)})\\
               &=\sum_{(a)}S_2(a_{(1)})a_{(2)}S_1(a_{(3)})\\
               &=\sum_{(a)}S_2(a_{(1)})a_{(2)}S_2(a_{(3)})\\
					&=S_2(a)
\endalign$$
for all $a$. This completes the proof.
\hfill$\square$
\einspr

In the proof above, one has to cover the formulas properly and one really has to  use that the equality $S_1=S_2$ means that $c(S_1(a)b)=(cS_2(a))b$ for all $a,b,c\in A$.
\snl
The result in the previous proposition is rather remarkable and turns out to be the key to the final definition. We explain this in the following remark. We refer to our forthcoming paper on the subject '{\it Weak Multiplier Hopf algebras I. The main theory}' ([VD-W2]) where proofs of the claims in the remark are given.

\iinspr{3.18} Remark \rm
i) First observe that the formulas in the proposition determine $F_1$ and $F_2$. This can be seen as follows. Assume e.g.\ that $a\in A$ and that $E(a\ot 1)=0$. This implies that $\Delta(b)(a\ot c)=0$ for all $b,c\in A$. Because the coproduct is assumed to be full, it follows that $da=0$ for all $d\in A$. Then $a=0$ because the product is assumed to be non-degenerate. From this it follows that $F_1$ is completely determined by the formula $E_{13}(F_1\ot 1)=E_{13}(1\ot E)$. Similarly we have that $a=0$ if $a\in A$ and $(1\ot a)E=0$. This will imply that $F_2$ is determined by the formula $(1\ot F_2)E_{13}  =(E\ot 1)E_{13}$.
\snl
ii) In fact, one has more. These two formulas can be used to (essentially) {\it define} the multipliers $F_1$ and $F_2$.
\snl
iii) Moreover, it is then relatively easy to obtain, using the observation in i), that these elements $F_1$ and $F_2$ have to satisfy the formulas in Proposition 3.15. The formulas will be a consequence of the given equations for $E$ (as in Proposition 3.9 and Assumption 3.10).
\snl
iv) It is also possible to show that $F_1$ and $F_2$ are idempotents (in the appropriate algebras) and that $1-F_1$ and $1-F_2$ are projection maps with range {\it in} the kernels of $T_1$ and $T_2$ respectively. If now it is assumed that they project {\it onto} these kernels, then the generalized inverses $R_1$ and $R_2$ are determined. They will automatically satisfy the Conditions 2.2 and give associated antipodes $S_1$ and $S_2$.
\snl
v) In Proposition 3.17 we obtained that the equality $S_1=S_2$ gave rise to the defining formulas for $F_1$ and $F_2$.  It turns out that if we obtain $S_1$ and $S_2$ as above, we can actually show that $S_1=S_2$. This is not so remarkable. However, it also follows that the antipode $S$ we get in this way, will be an anti-algebra and an anti-coalgebra map. This will not be shown here, but in the [VD-W2]. 
\hfill$\square$\einspr

This all means that, given the conditions on the ranges of $T_1$ and $T_2$ in terms of the idempotent $E$ (with the right properties), the requirement $S_1=S_2$ determines the kernels of these maps and a good choice of projection maps on these kernels, giving an antipode $S$ with the right properties. We believe that this is a quite {\it remarkable result}.
\nl
In the next section, where we treat the regular case, we will see how all these formulas behave with respect to the involutive structure in the case of a $^*$-algebra. 
\nl
We finish the section, as promised, with a look at the examples.

\iinspr{3.19} Example \rm
i) For the case $A=K(G)$, we have seen in Example 3.6.i that $E$ is the element in $M(A\ot A)$ given by the function that is $1$ on pairs $(p,q)$ for which $s(p)=t(q)$. Then $F_1$ is given by the function that is $1$ on pairs $(p,q)$ for which $s(p)=s(q)$ and $F_2$ is given by the function that is $1$ on pairs $(p,q)$ for which $t(p)=t(q)$. We leave it to the reader as an exercise to verify the various results and formulas involving these idempotents. However, let us just look at the first formula in Proposition 3.17. We get in this case
$(E_{13}(F_1\ot 1))(p,q,v)=1$ if and only if $s(p)=t(v)$ and $s(p)=s(q)$. On the other hand we get
$(E_{13}(1\ot E))(p,q,v)=1$ if and only if $s(p)=t(v)$ and $s(q)=t(v)$. These conditions are the same.
\snl
ii) For the case $A=\Bbb CG$, we have seen in Example 3.6.ii that $E=\sum_e \lambda_e\ot \lambda_e$ where the sum is taken over all the units of the groupoid. Because now $S(\lambda_e)=\lambda_e$ for every unit, we get that $F_1$ and $F_2$ are given by the same expression. Again we leave it to the reader as an exercise to verify the various results and formulas. 
\hfill$\square$\einspr

\iinspr{3.20} Example \rm
Consider now the case of a weak Hopf algebra. We know that $E=\Delta(1)$ and $F_1=(\iota\ot S)\Delta(1)$ and $F_2=(S\ot \iota)\Delta(1)$ where $S$ is the antipode. Let us again consider only the first formula of Proposition 3.17. In Proposition 2.3.4 of [N-V2], we find that $E(x\ot 1)=E(1\ot S(x))$ when $x\in A_s$. Remember that $A_s$ is the left leg of $E$. Then we find indeed
$$E_{13}(F_1\ot 1)=(\iota\ot S\ot\iota)(E_{13}(E\ot 1))=E_{13}(1\ot E)$$
because $\sigma (S\ot S)E=E$.
\hfill$\square$\einspr

\nl



\bf 4. The definition of a (regular) weak multiplier Hopf algebra \rm
\nl
We are now ready to give (a first version of) the definition of a weak multiplier Hopf algebra. The approach is by generalizing the original definition of a multiplier Hopf algebra. Further in this section, we will also consider the regular case.  In our second work on the subject [VD-W2], we will  consider other possible (equivalent) definitions for a (regular) weak multiplier Hopf algebra. 
\nl
\it The definition of a weak multiplier Hopf algebra \rm
\nl
We assume that $A$ is a non-degenerate idempotent algebra with a full coproduct $\Delta$ such that there exists a counit $\varepsilon$. We know that the counit is uniquely determined because the coproduct is assumed to be full (see Proposition 1.12 in Section 1). 
\snl
We consider the maps $T_1$ and $T_2$ from $A\ot A$  to itself as defined in Section 1 by
$$T_1(a\ot b)=\Delta(a)(1\ot b)\qquad\text{and}\qquad T_2(a\ot b)=(a\ot 1)\Delta(b)$$
for $a,b\in A$. As before, we denote the ranges of $T_1$ and $T_2$ by $\text{Ran}(T_1)$ and $\text{Ran}(T_2)$ respectively and the kernels by $\text{Ker}(T_1)$ and $\text{Ker}(T_2)$. 
\snl
In the case of a multiplier Hopf algebra these maps are assumed to be bijective. This is no longer the case here. The following is {\it a definition} (or perhaps rather a characterization) of a weak multiplier Hopf algebra {\it in terms of the ranges and kernels of these canonical maps}. We will consider (a modified version of) this definition again in [VD-W2] where we take it as the starting point of the development of weak multiplier Hopf algebras.

\inspr{4.1} Definition \rm Let $(A,\Delta)$ be a pair of an algebra with a coproduct as above. Assume that there is an idempotent multiplier $E$ in $M(A\ot A)$ such that  
$$\text{Ran}(T_1)=E(A\ot A)\qquad\text{and}\qquad \text{Ran}(T_2)=(A\ot A)E \tag"(4.1)"$$  
and that it satisfies 
$$(\iota\ot\Delta)(E)=(\Delta\ot\iota)(E)=(1\ot E)(E\ot 1)=(E\ot 1)(1\ot E).\tag"(4.2)"$$
Let $F_1$ be a right multiplier of $A\ot A^{\text{op}}$ and $F_2$ a left multiplier of $A^{\text{op}}\ot A$ such that 
$$ E_{13}(F_1\ot 1) = E_{13}(1\ot E) \qquad\text{and}\qquad (1\ot F_2)E_{13} =(E\ot 1)E_{13}\tag"(4.3)"$$
and that  
$$\text{Ker}(T_1)=(A\ot 1)(1-F_1)(1\ot A) \qquad\text{and}\qquad \text{Ker}(T_2)=(A\ot 1)(1-F_2)(1\ot A).\tag"(4.4)"$$
Then we call $(A,\Delta)$ a {\it weak multiplier Hopf algebra}.
\hfill$\square$
\einspr

As this definition is very important in this paper, and not so obvious, we now make a couple of important remarks.

\inspr{4.2} Remark \rm
i) We have seen that condition (4.1) uniquely determines $E$ (see Lemma 3.3 and Proposition 3.5 in Section 3) and that it makes the extension of $\Delta$ as well as of $\iota\ot\Delta$ and $\Delta\ot\iota$ to the multiplier algebras possible (see Appendix A). So condition (4.2) makes sense in $M(A\ot A\ot A)$. We also know that condition (4.2) encodes the commutation rules of the coproduct with the maps $T_1$ and $T_2$, as well as with the maps $T_1R_1$ and $T_2R_2$ (see Proposition 3.9 and 3.11 in Section 3). 
\snl
ii) We have seen that the formulas (4.3) make sense and characterize these multipliers $F_1$ and $F_2$ (see Remark 3.18.i) so that they are completely determined by $E$ (and hence by the coproduct itself). In fact, roughly speaking, the formulas (4.3) can be used to define the multipliers $F_1$ and $F_2$ (see [VD-W2]). 
We also know that $F_1$ and $F_2$ automatically satisfy the necessary commutation rules with $\Delta$. 
\snl
iii) Remark that the first formula in condition (4.4) is equivalent with the property
$$T_1(\sum_i a_i\ot b_i)=0 \qquad \Longleftrightarrow \qquad  \sum_i (a_i\ot 1)F_1(1\ot b_i)=0$$
for elements $a_i, b_i$ in $A$. A similar result is true for the second formula.
\snl 
iv) Because $E$, $F_1$ and $F_2$ are uniquely determined by the coproduct itself, we do not need to include these objects in the notation for a weak multiplier Hopf algebra. The same is true for the counit $\varepsilon$ and the antipode $S$.
\snl
v) Remark that if $E=1$, that is when the maps $T_1$ and $T_2$ are assumed to be surjective, then we must have by the assumptions that $F_1=1$ as well as $F_2=1$ so that actually, we have a multiplier Hopf algebra.
\snl
vi) Finally, remark that formulas (4.2) and (4.3) only involve the legs of $E$, $F_1$ and $F_2$. This would allow us to give (still) another approach to the theory.
\hfill$\square$
\einspr

As we mentioned already, in [VD-W2] we will start from a modified version of this definition and build the theory. The formulas in (4.4) together with (4.1) determine generalized inverses $R_1$ and $R_2$ of $T_1$ and $T_2$. The formulas in (4.2) give the necessary commutation rules with the coproduct for the ranges and those in (4.3) give, together with the first ones in (4.2) also the necessary commutation rules with the coproduct for the kernels. Then we can obtain the associated antipodes $S_1$ and $S_2$. They will satisfy $S_1=S_2$. The antipode $S$ thus defined is both an anti-homomorphism and it flips the coproduct. 
\snl
We refer to [VD-W2] for the proofs (and more details about this reasoning).
\nl
We now look at the {\it basic examples}.  We will collect the results and show that for all cases we considered, we do get a weak multiplier Hopf algebra in the sense of Definition 4.1 above. Most of the properties we need have been considered already in the previous sections.
\snl
First we consider the case of functions on a groupoid.  

\inspr{4.3} Proposition \rm
Let $G$ be a groupoid. Let $A$ be the algebra of complex functions with finite support on $G$. Define a coproduct $\Delta$ on $A$ by $\Delta(f)(p,q)=f(pq)$ if $p,q\in G$ and if $pq$ is defined. Otherwise we let $\Delta(f)(p,q)=0$. Then $(A,\Delta)$ is a weak multiplier Hopf algebra (in the sense of Definition 4.1). The antipode is given by $S(f)(p)=f(p^{-1})$ for all $f\in A$ and $p\in G$. 

\snl\bf Proof\rm:
We have already seen in Section 1 that $\Delta$, as defined above on $A$, is a full coproduct with a counit on the algebra $A$ (see Proposition 1.6 and Example 1.13.i).
\snl
In Section 3, we have defined $E$ as a function on $G\times G$ by $E(p,q)=1$ if $p,q\in G$ and if $s(p)=t(q)$. Otherwise we have put $E(p,q)=0$. It is an idempotent in $M(A\ot A)$ satisfying the required conditions (4.1) and the formulas in (4.2) in Definition 4.1. See Example 3.6.i and Example 3.19.i.
\snl
Then we define $F_1$ and $F_2$ in $M(A\ot A)$. Remark that the algebra is abelian so that $A^{\text{op}}$ and $A$ are the same. For $p,q\in G$ we put
$F_1(p,q)=1$ if $s(p)=s(q)$ and $F_1(p,q)=0$ otherwise. Similarly, we put
$F_2(p,q)=1$ if $t(p)=t(q)$ and $F_2(p,q)=0$ otherwise. As we mentioned already in Example 3.19.i, a simple argument will show that these idempotents satisfy the equations (4.3) in Definition 4.1. We have e.g.\ that
$(E_{13}(F_1\ot 1))(u,v,w)=1$ if and only if $s(u)=t(w)$ and $s(u)=s(v)$. Similarly, we have $
(E_{13}(1\ot E))(u,v,w)=1$ if and only if $s(u)=t(w)$ and $s(v)=t(w)$. These two conditions are the same.
\snl
One can also check the conditions (4.4) in Definition 4.1. Suppose e.g.\newline   $\sum_i\Delta(a_i)(1\ot b_i)=0$ for a finite number of elements $a_i,b_i\in A$. Then, if $p,q\in G$ and if $pq$ is defined, we find
$$\sum_i(\Delta(a_i)(1\ot b_i))(p,q)= \sum_i a_i(pq)b_i(q).$$
This means that  $\sum_i a_i(u)b_i(q)=0$ whenever $s(u)=s(q)$. On the other hand
$$\sum_i((a_i\ot 1)F_1(1\ot b_i))(u,q)=\sum_i((a_i\ot b_i))(u,q)$$ for all $u,q\in G$ such that $s(u)=s(q)$. From all this, we see very easily that $\sum_i\Delta(a_i)(1\ot b_i)=0$ if and only if $\sum_i(a_i\ot 1)F_1(1\ot b_i)=0$. This will give the first equation in (4.4) of Definition 4.1. The argument for the other one is very similar.
\snl
This shows that we  have a weak multiplier Hopf algebra as in Definition 4.1. 
\snl
Furthermore, it is easy to show that the associated antipode is given by $(S(f))(p)=f(p^{-1})$ as expected. And then we see that $F_1=(\iota\ot S)E$ and $F_2=(S\ot \iota)E$ (because $t(p^{-1})=s(p)$ for all $p\in G$). See also Section 2, Example 2.14.i.
\hfill$\square$\einspr

The dual case is very similar.

\inspr{4.4} Proposition \rm 
Let $G$ be a groupoid. Let $B$ be the groupoid algebra $\Bbb CG$ and define $\Delta$ on $B$ by $\Delta(\lambda_p)=\lambda_p\ot \lambda_p$ where $p\mapsto \lambda_p$ is the imbedding of $G$ in the groupoid algebra $\Bbb CG$. Then $(B,\Delta)$ is a weak multiplier Hopf algebra. The antipode $S$ is given by the formula $S(\lambda_p)=\lambda_{p^{-1}}$ for all $p\in G$.

\snl\bf Proof\rm: 
We have shown already in Proposition 1.7 that $\Delta$ is a regular and full coproduct on $B$. There is also a counit.
\snl 
In Section 3, we have obtained the multiplier $E$ also for this case. It is given as $\sum_e \lambda_e\ot \lambda_e$ where the sum is taken over all units. In Example 3.6.ii and Example 3.19.ii it is shown that this multiplier satisfies the conditions (4.1) and the formulas (4.2) in Definition 4.1. We have defined $F_1$ and $F_2$ both using the same formula as for $E$ and we have mentioned in Example 3.19.ii that we get multipliers as in Definition 4.1, satisfying the formulas (4.3).
\snl
Again, essentially, the only thing left to check is condition (4.4). 
\snl
Assume e.g.\ that $\sum_i \Delta(a_i)(1\ot b_i)=0$ with a finite number of elements $(a_i)$ and $(b_i)$ in $B$. Write 
$$a_i=\sum_p a_i(p)\lambda_p
	\qquad\text{and}\qquad\quad
		b_i=\sum_q b_i(q)\lambda_q$$
for all $i$.
Then 
$$\sum_i \Delta(a_i)(1\ot b_i)=\sum_{i,p,q}a_i(p)b_i(q)\lambda_p\ot \lambda_{pq}$$
where the sum is only taken over those pairs $(p,q)$ for which $pq$ is defined. If this is equal to $0$, then we must have that $\sum_i a_i(p)b_i(q)=0$ for all pairs $(p,q)$ with $s(p)=t(q)$.
On the other hand, we have
$$\align \sum_i (a_i\ot 1)F_1(1\ot b_i)
	&=\sum_{i,p,q,e} a_i(p)b_i(q)\lambda_{pe}\ot \lambda_{eq}\\
	&=\sum_{i,p,q} a_i(p)b_i(q)\lambda_{p}\ot \lambda_{q}
\endalign$$
where the last sum is only taken over those pairs for which $s(p)=t(q)$. Again we conclude easily that $\sum_i \Delta(a_i)(1\ot b_i)=0$ if and only if $\sum_i (a_i\ot 1)F_1(1\ot b_i)=0$.
\snl
Therefore, we have shown that the pair $(B,\Delta)$ is also a weak multiplier Hopf algebra.
\snl
In Example 2.14.ii we have seen that the antipode is given by $S(\lambda_p)=\lambda_{p^{-1}}$ for all $p\in G$.
\hfill$\square$\einspr

Later in this section, we consider the case of weak Hopf algebras.
\nl
\it Regular weak multiplier Hopf algebras \rm
\nl
As we have noticed in previous discussions, there are some peculiarities in the non-regular case that are badly understood. They disappear when the weak multiplier Hopf algebra is assumed to be regular as we will see in what follows. 
\snl
Let us now assume that we have a weak multiplier Hopf algebra (as in Definition 4.1), let $R_1$ and $R_2$ be the unique generalized inverses of $T_1$ and $T_2$ determined by $E$, $F_1$ and $F_2$ and let $S$ be the associated antipode. We know that $S:A \to M(A)$ and that it is an anti-algebra map and an anti-coalgebra map (proven in [VD-W2]). 
\snl 
We arrive at the following definition of a  regular weak multiplier Hopf algebra.

\inspr{4.5} Definition \rm
Let $(A,\Delta)$ be a weak multiplier Hopf algebra. We call it {\it regular} if the antipode $S$ maps $A$ to itself and if it is bijective.
\hfill$\square$\einspr

This is indeed what we expect for a regular weak multiplier Hopf algebra, given the notion of regularity for ordinary multiplier Hopf algebras. Later in this section, we will find necessary and sufficient conditions for regularity in terms of the coproduct. Then we will also consider the case of weak multiplier Hopf $^*$-algebras and prove that they are automatically regular.
\snl
Of course, also all the examples we have considered up to now are regular because in the groupoid case, the antipode $S$ satisfies $S^2=\iota$ whereas in the case of a finite-dimensional weak Hopf algebra, the antipode is proven to be bijective (see e.g.\ Proposition 2.3.1 in [N-V2] and Proposition 2.10 in [B-N-S]).
\nl
In the previous two sections, we have already considered the assumption that $S$ maps $A$ bijectively to itself. See e.g.\ Proposition 2.10 in Section 2 and also several remarks in Section 3. This was done for motivational reasons. The difference here is that now we have it as a real assumption and the consequences are now genuine properties for a regular weak multiplier Hopf algebra.
\snl
In the first place, we see what can be concluded about the multipliers $E$, $F_1$ and $F_2$.
\snl
We first have the following expected property.

\inspr{4.6} Proposition \rm 
Let $(A,\Delta)$ be a regular weak multiplier Hopf algebra.
Then we have $(S\ot S)E=\sigma E$ where as before $\sigma$ is the flip on $A\ot A$, extended to $M(A\ot A)$. Similarly $(S\ot S)F_2=\sigma F_1$.
\hfill$\square$\einspr

The proof is rather straightforward and we will give details in [VD-W2]. 
\snl
The formulas in the formulation of the previous proposition can be interpreted using the way they are proven above. It is also possible to extend $S\ot S$ to the multiplier $M(A\ot A)$ using the known techniques. This is then used for the first formula. Similarly, $S\ot S$ can be extended  to $M(A^{\text{op}}\ot A)$ and this is used for the second formula.
\snl
The above two formulas are basically natural consequences of the simple fact that the antipode converts the map $T_1$ to $T_2$. However, the antipode can also be used to get a relation between the maps $R_1$ and $R_2$ on the one hand and the maps $T_3$ and $T_4$ on the other hand. We have considered this relation already in Propositions 2.11 and 2.12 in Section 2 and in Remark 3.13 in the previous section. Indeed we have 
$$ R_1 (\iota \ot S) = (\iota\ot S) T_3 
	\qquad\text{and}\qquad\quad
		R_2 (S\ot \iota)  = (S\ot \iota) T_4.
$$
Therefore, as we also have observed already, we expect the following formulas given $F_1$ and $F_2$ in terms of $E$. We will only formulate the result. For the proof, we refer again to [VD-W2]. 

\inspr{4.7} Proposition \rm
Let $(A,\Delta)$ be a regular weak multiplier Hopf algebra. Then
$$F_1=(\iota\ot S)E 
	\qquad\text{and}\qquad\quad
		F_2=(S\ot \iota)E.$$
\vskip -0.5 cm \hfill$\square$\einspr

Remark that the formulas in Proposition 4.7 are completely in accordance with the ones in Proposition 4.6.
\snl
For a correct interpretation of e.g.\ the formula $F_1=(\iota\ot S)E$, one can write it as
$$(a\ot 1)F_1(1\ot S(b))=(\iota\ot S)((a\ot b)E)$$
for all $a,b\in A$. Another possibility interpretation is by extending the map $\iota\ot S$ to $M(A\ot A)$ first, but that is essentially the same story.
\snl
Similarly for the other formula.
\snl
From the formulas in Proposition 4.6 and 4.7, we can now completely formulate the data for the two other canonical maps $T_3$ and $T_4$. First we have the following.

\inspr{4.8} Proposition \rm
Let $(A,\Delta)$ be a regular weak multiplier Hopf algebra. Then the coproduct is regular (as in Definition 1.3) and we can find generalized inverses $R_3$ and $R_4$ of the maps $T_3$ and $T_4$ respectively given by the formulas
$$\align 
	R_3 &= (\iota \ot S^{-1}) T_1 (\iota\ot S) \\
	R_4 &= (S^{-1}\ot \iota)  T_2  (S\ot \iota).
\endalign$$
The associated (inverse) antipodes $S_3$ and $S_4$ exist and satisfy
$S_3=S_4=S^{-1}$.

\snl \bf Proof\rm:
The proof is rather straightforward. If e.g.\ we apply $R_3$ to $a\ot b$ with $a,b\in A$ and use the formula above, we get
$$R_3(a\ot b)=\sum_{(a)}a_{(1)} \ot bS^{-1}(a_{(2)})$$
and from the definition of $S_3$ (see Section 2), we find $S_3=S^{-1}$. Similarly we get $S_4=S^{-1}$ from the formula above defining $R_4$ and the definition of $S_4$ as given in Section 2.
\hfill $\square$\einspr

For the associated projection maps, we find the following formulas.

\inspr{4.9} Proposition \rm With these choices of $R_3$ and $R_4$ we find that, 
$$T_3R_3 (a\ot b) = (a\ot b)E 
	\qquad\quad\text{and}\qquad\quad
	T_4R_4 (a\ot b) = E(a\ot b)
$$
for all $a,b\in A$.
We also have
$$ R_3T_3 (a\ot b) = (1\ot b)F_3 (a\ot 1)
	\qquad\quad\text{and}\qquad\quad
		R_4T_4 (a\ot b) = (1\ot b)F_4 (a\ot 1)
$$
where $F_3$ is the left multiplier of $A\ot A^{\text{op}}$ given by $F_3=(\iota\ot S^{-1})E$ and $F_4$ is the right multiplier of $A^{\text{op}}\ot A$ given by $F_4=(S^{-1}\ot \iota)E$.  

\snl\bf Proof\rm:
The proof of these four formulas follows immediately from the given relations between the pair $(T_3, T_4)$ and the pair $(R_1,R_2)$ and the associated relation between the pairs $(T_1,T_2)$ and $(R_3,R_4)$.
\hfill $\square$ \einspr  

Remark that this is what we expected already. Remember that in Section 3, we did use the formulas above for $T_3R_3$ and $T_4R_4$ as a starting point for further investigations (see 
again Remark 3.13). This assumption led us there to the formulas for $F_1$ and $F_2$ given in Proposition 4.7 above.

Combining the various results we obtained already in this section,  we find the following collection of formulas.
\snl
First, we can express the idempotents $F_1, F_2, F_3$ and $F_4$ all in terms of $E$:
$$\align F_1&=(\iota\ot S)E  \qquad\quad\text{and}\qquad\quad F_3=(\iota\ot S^{-1})E \tag"(4.5)"\\
	F_2&=(S\ot \iota)E \qquad\quad\text{and}\qquad\quad F_4=(S^{-1}\ot \iota)E.\tag"(4.6)"
\endalign$$
On the other hand, we also have the following four formulas:
$$\align E_{13}(F_1\ot 1) &= E_{13}(1\ot E) 
	\qquad\quad\text{and}\qquad\quad 
		(F_3\ot 1)E_{13}=(1\ot E)E_{13} \tag"(4.7)"\\
	(1\ot F_2)E_{13} &=(E\ot 1)E_{13}
	\qquad\quad\text{and}\qquad\quad 
	 	E_{13}(1\ot F_4)=E_{13}(E\ot 1).\tag"(4.8)"
\endalign$$
The formulas with $F_1$ and $F_2$ are part of the definition and the formulas with $F_3$ and $F_4$ follow by applying the antipode and using the various relations above.
\snl

Before we deduce from this an equivalent definition of a regular weak multiplier Hopf algebra, we first want to make some important remarks about the symmetry we discover above.

\iinspr{4.10} Remark \rm
i) Assume that $(A,\Delta)$ is a regular weak multiplier Hopf algebra. Consider a new pair $(A^{\text{op}},\Delta)$ where $A^{\text{op}}$ is the algebra $A$ but with the opposite product and with the same coproduct. With this procedure, the original maps $T_3,T_4$ are replaced by $T_1,T_2$ for the new pair. We see first that the multiplier $E$ does not change. Further, if we consider the formulas (4.7) and (4.8) above, we see that passing from $A$ to $A^{\text{op}}$ will interchange the two formulas in (4.7) as well as the two formulas in (4.8), precisely as expected.
\snl
On the other hand, if we consider the formulas (4.5) and (4.6), we see that we have to replace $S$ by $S^{-1}$. Also this is expected. 
\snl
ii) Now, consider the transition from $(A,\Delta)$ to $(A,\Delta^{\text{cop}})$ instead. This is slightly more complicated. The original maps $T_3,T_4$ are now replaced by $\sigma T_2\sigma, \sigma T_1\sigma$. This implies that $E$ will become $\sigma E$ for the new pair. Further, consider the formulas (4.7) and (4.8). If we e.g.\ apply $\sigma_{13}$ to the equality $E_{13}(1\ot F_4)=E_{13}(E\ot 1)$ we find $(\sigma E)_{13}((\sigma F_4)\ot 1) = (\sigma E)_{13}(1\ot (\sigma E))$. This is precisely as it should because $E$ is replaced by $\sigma E$ and $F_4$ by $\sigma F_1$. Similarly, for the formula with $F_3$.
\snl
Again, if we look at the formulas (4.5) and (4.6), we see that also here, we have to replace $S$ by $S^{-1}$.
\hfill$\square$\einspr

One can also verify the transition from $(A,\Delta)$ to  $(A^{\text{op}},\Delta^{\text{cop}})$ as we did above. Then the pair $(T_1,T_2)$ is replaced by the pair $(T_2, T_1)$. This will not give anything new.
\snl
We see from all these observations that the pieces of the puzzle fit very nicely together. 
\nl
It also suggests the following {\it equivalent characterization} of regular weak multiplier Hopf algebras.

\iinspr{4.11} Proposition  \rm
Let $(A,\Delta)$ be a weak multiplier Hopf algebra. Then it is regular if and only if also $(A^{\text{op}},\Delta)$ (or equivalently $(A,\Delta^{\text{cop}})$) is a weak multiplier Hopf algebra.

\hfill$\square$\einspr

One direction is immediately clear from the results that we have obtained. The converse will be shown in [VD-W2]. Indeed, we will prove that the antipode is a bijection from $A$ to itself if also $(A^{\text{op}},\Delta)$ (or equivalently $(A,\Delta^{\text{cop}})$) is a weak multiplier Hopf algebra.
\snl
As an immediate consequence, we will have that if $(A,\Delta)$ is a  weak multiplier Hopf algebra with either $A$ is abelian of $\Delta$ coabelian, then it is regular. In these two cases, we get $S=S^{-1}$.
\snl
Another consequence is that the two weak multiplier Hopf algebras, associated to a group\-oid are automatically regular. Furthermore, also a finite-dimensional weak Hopf algebra is regular because the antipode is bijective in that case.
\nl
In the next proposition, we formulate various results about the relation of weak Hopf algebras with weak multiplier Hopf algebras. We will not be able to prove all the statements, but for the missing arguments, we refer to [VD-W2].

\iinspr{4.12} Proposition \rm i) Let $(A,\Delta)$ be a weak Hopf algebra. Then it is a weak multiplier Hopf algebra.\newline
ii) Conversely, if $(A,\Delta)$ is a {\it regular} weak multiplier Hopf algebra and if the underlying algebra $A$ is unital, it is a weak Hopf algebra.

\snl\bf Proof \rm (sketch):
i) Assume first that $(A,\Delta)$ is a weak Hopf algebra. There is a counit by assumption and because $\Delta$ maps $A$ to $A\ot A$, the coproduct is automatically full.
\snl
With $E=\Delta(1)$, we have seen already in Section 3 that the condition (4.1) and (4.2) of Definition 4.1 are fulfilled (see Example 3.7).
\snl
If we let $F_1=(\iota\ot S)E$ and $F_2=(S\ot \iota)E$ we get idempotent elements in $A\ot A^{\text{op}}$ and $A^{\text{op}}\ot A$ respectively because $S$ is an anti-homomorphism. In Example 3.20, we have shown how (4.3) of Definition 4.1 follows.
\snl
Finally, we show that also (4.4) in Definition 4.1 is satisfied. So, again assume that $\sum_i \Delta(a_i)(1\ot b_i)=0$ with a finite number of elements $(a_i)$ and $(b_i)$ in $A$. If we apply $R_1$ with the antipode $S$, we find $\sum_i (a_i\ot 1)F_1(1\ot b_i)=0$ because we know that $R_1T_1(a\ot b)=(a\ot 1)F_1(1\ot b)$ for all $a,b$ as we have seen in Section 3 (see Example 3.20).
\snl
In [VD-W2] we give a more elegant argument using an alternative definition for a weak multiplier Hopf algebra (see Theorem 2.9 and Proposition 2.10 in [VD-W2]).
\snl
ii) Conversely, consider a regular weak multiplier Hopf algebra $(A,\Delta)$ and assume that $A$ has an identity. It is shown that there exists an antipode and it is easy to verify that it satisfies the required conditions. What still has to be shown is the  so-called {\it weak multiplicativity} of the counit (see Definition 2.1 in [B-N-S]). The argument goes as follows.
\snl
We start with 
$$(1\ot a)\Delta(b)(c\ot 1)=\sum_{(c)}(1\ot a)\Delta(bc_{(1)}(1\ot S(c_{(2)})),$$ 
true for all $a,b,c\in A$. If we apply $\varepsilon\ot\varepsilon$ we find
$$(\varepsilon\ot\varepsilon)((1\ot a)\Delta(b)(c\ot 1))=\sum_{(c)}\varepsilon(abc_{(1)}S(c_{(2)})).$$  
If we apply this with $a=1$ we find 
$$\varepsilon(bc)=(\varepsilon\ot\varepsilon)\Delta(b)(c\ot 1)=\sum_{(c)}\varepsilon(bc_{(1)}S(c_{(2)}))$$
and if we use this formula with $b$ replaced by $ab$ in the previous formula, we find
$$(\varepsilon\ot\varepsilon)((1\ot a)\Delta(b)(c\ot 1))=\varepsilon(abc)$$
for all $a,b,c\in A$. This gives one of the properties of the counit we need.
\snl
Now we use the assumption that $A$ is regular. Then we can apply the previous result for $(A,\Delta^{\text{cop}})$ and since the counit is the same, we find the other formula
$$(\varepsilon\ot\varepsilon)((a\ot 1)\Delta(b)(1\ot c))=\varepsilon(abc)$$
for all $a,b,c\in A$. We also find this formula from the other one, applied to  $(A^{\text{op}},\Delta)$.

\hfill$\square$\einspr

When in this paper, we refer to a weak Hopf algebra, we have the references [B-N-S] and [N-V2] in mind. In these papers, only the finite-dimensional case is considered. In the finite-dimensional case, the antipode is proven to be invertible (see e.g.\ Theorem 2.10 in [B-N-S]). This implies that finite-dimensional weak Hopf algebras are automatically regular. Conversely, if we have a regular weak multiplier Hopf algebra with a finite-dimensional underlying algebra, it has to be a weak Hopf algebra. This statement follows from item ii) in the previous proposition. One has to argue that the underlying algebra of a finite-dimensional regular weak multiplier Hopf algebra has to be unital, but this follows from the existence of local units (see Proposition 4.9 in [VD-W2]).
\snl
A second remark is the following. The statement in the previous proposition is not completely symmetric. We seem to need regularity to go back. Indeed, when we look closer to the various arguments, we see that the weak multiplicativity of the counit is not used to prove that any weak Hopf algebra is a weak multiplier Hopf algebra. On the other hand, conversely, without regularity we can only prove one weak multiplicativity formula for the counit. In order to prove the second one, we  need regularity. We refer to the results and the discussions about this in Section 4 of [VD-W2]. 
\nl
We now  finish this section by the involutive case as this is directly linked with the regularity.
\snl
If $(A,\Delta)$ is a weak multiplier Hopf algebra and $A$ is a $^*$-algebra and $\Delta$ a $^*$-homomorphism, we will have $E=E^*$ as we have seen in Section 3. Then the following definition makes sense.

\iinspr{4.13} Definition \rm
Let $(A,\Delta)$ be a weak multiplier Hopf algebra and assume that $A$ is a $^*$-algebra and $\Delta$ a $^*$-homomorphism. Then we call $(A,\Delta)$ a {\it weak multiplier Hopf $^*$-algebra}.
\hfill$\square$\einspr

Also the next proposition is expected.

\iinspr{4.14} Proposition \rm
If $(A,\Delta)$ is a weak multiplier Hopf $^*$-algebra, then it is regular. The antipode satisfies $S(S(a)^*)^*=a$ for all $a\in A$. And not only do we have $E^*=E$ but also 
$$F_1^*=F_3 \qquad\quad\text{and}\qquad\quad F_2^*=F_4.$$

\snl\bf Proof\rm:
We know that the coproduct is regular and that 
$$T_3(a^*\ot b^*)=T_1(a\ot b)^*
     \qquad\quad\text{and}\qquad\quad
       T_4(a^*\ot b^*)=T_2(a\ot b)^*. 
$$
From this it easily follows that $(A^{\text{op}},\Delta)$ is again a weak multiplier Hopf algebra and that the idempotent $E$ is the same for both $(A,\Delta)$ and $(A^{\text{op}},\Delta)$ (also because $E^*=E$). Therefore, also $(A,\Delta)$ it is a regular weak multiplier Hopf algebra. As we have seen already before, we will have that the antipode of $(A^{\text{op}},\Delta)$ is $S^{-1}$ and also given by $a\mapsto S(a^*)^*$ (cf.\ Proposition 2.9) in Section 2.  This will imply the property of $S$ as in the formulation of the proposition. Finally, the equalities $F_1^*=F_3$ and  $F_2^*=F_4$ follow from the formulas 
$$R_3(a^*\ot b^*)=R_1(a\ot b)^*
     \qquad\quad\text{and}\qquad\quad
       R_4(a^*\ot b^*)=R_2(a\ot b)^* 
$$
and the definitions of these idempotents.
\hfill$\square$
\einspr 

The weak multiplier Hopf algebras that are obtained from a groupoid as in Proposition 4.3 and Proposition 4.4 are weak multiplier Hopf $^*$-algebras as it is checked easily that the coproduct are $^*$-homomorphisms. 

\nl\nl



\bf 5. Conclusions and final remarks\rm
\nl
In this paper, we have step by step developed a possible definition of what we call a {\it weak multiplier Hopf algebra}. The result is found in Definition 4.1 in Section 4. We also have given the definition of a {\it regular} weak multiplier Hopf algebra (see Definition 4.5). 
\snl
We have given several arguments for the claim that these are good and natural definitions. We have always tried to have conditions that  naturally are self-dual in the situation of a dual pair. And of course, we have argued that the basic examples, coming from a groupoid, as well as the known case of a weak Hopf algebra, fit into our theory.
\snl
There are two basic aspects in this theory. First there is the antipode describing the generalized inverses $R_1$ and $R_2$ of the canonical maps. The antipode is uniquely determined by these inverses. The other aspect concerns the idempotents that determine the choice of these inverses. Properties of these idempotents take place on a different level. That is, roughly speaking, about the legs of the idempotent $E$ (playing the role of $\Delta(1)$). This second aspect is typical for the theory of weak (multiplier) Hopf algebras and it is not present in the case of (multiplier) Hopf algebras (as the canonical maps are then assumed to be bijective). We can say that the first aspect is treated in Section 2 of the paper, while the second one is investigated in Section 3.
\nl
Of course, a lot of work still has to be done. First we need to develop the theory from the definition (as in this paper, we rather have 'developed the definition from the theory'). This will be done in our second paper on the subject entitled {\it Weak multiplier Hopf algebras I. The main theory}, see reference [VD-W2]. In that paper, we also give the proofs of some of the results that we have stated in this paper without proof.
In [VD-W2], we will also start with the study of the {\it  source and target algebras} $A_s$ and $A_t$, as well as the {\it source and target maps} $\varepsilon_s:A \to A_s$ and $\varepsilon_t:A\to A_t$. These algebras will be subalgebras of the multiplier algebra $M(A)$ (and not of $A$ in general). 
\snl
In a third paper on the subject, entitled {\it Weak multiplier Hopf algebras II. The source and target algebras}, we will investigate the objects further. And we will use the results to look at examples (see [VD-W3]).
\snl
Finally, in still another paper {\it Weak multiplier Hopf algebras III. Integrals and duality}, we will study integrals on weak multiplier Hopf algebras and show that the reduced dual of $A$ is again a weak multiplier Hopf algebra (see [VD-W3]).
\nl
There are also various aspects of the theory that are not very well understood. 
\snl
First there are the different requirements for the underlying algebra. It is quite natural to assume that the algebra is non-degenerate. But what about the assumption that it is idempotent? We know that in the case of multiplier Hopf algebras, this is a property that can be proven from the axioms. And only in the regular case, we were able to show that the algebra has local units. Again in the case of multiplier Hopf algebras, this can be shown also in the non-regular case.
\snl
Secondly, there are some questions about the coproduct. We need that it is full in order to get uniqueness of the counit. This is not necessary in the case of a multiplier Hopf algebra because then the  existence of a counit implies fullness of the coproduct. It is also not clear what the conditions on the coproduct mean on the dual in the case of a dual pair.
\snl
Finally, there several peculiarities in the non-regular case that are badly understood. We have e.g.\ no example of a coproduct that is not a regular coproduct. We refer to Section 5 in [VD-W2] for more related comments. 
\nl\nl



\bf Appendix A. Extension of the coproduct to the multiplier algebra \rm
\nl
In this appendix {\it we assume} that $(A,\Delta)$ is a pair of a non-degenerate algebra $A$ with a coproduct $\Delta$ as in Definition 1.1 of this paper. We also assume that the algebra is idempotent, i.e.\ that $A^2=A$. This last condition is assumed to be satisfied in this paper (see a remark in Section 1).
\snl
If the coproduct is non-degenerate (as it is the case for multiplier Hopf algebras), there is a unique extension of $\Delta$ to the multiplier algebra $M(A)$ of $A$. The extension is unital. In the situation  we consider in this paper however, we cannot assume that the coproduct is non-degenerate and so we can not apply this result. In this appendix, we will find a {\it workaround} for this problem and we will obtain a generalization of this result that is useful for the study of weak multiplier Hopf algebras.
\snl
We need a condition on the coproduct that is {\it weaker} than non-degeneracy.

\ainspr{A.1} Assumption \rm 
We  assume that there is an idempotent $E$ in $M(A\ot A)$ such that
$$ E(A\ot A)=\Delta(A)(A\ot A) \qquad\qquad\text{and}\qquad\qquad
         (A\ot A)E=(A\ot A)\Delta(A).$$
\einspr

We have seen in Lemma 3.3 that if such an idempotent exists, then it is unique. So the assumption is in fact a {\it condition} on the coproduct. Remark that as in Proposition 3.5, also here we will get that $E\Delta(a)=\Delta(a)$ and $\Delta(a)E=\Delta(a)$ for all $a\in A$ and that $E$ is the smallest idempotent in $M(A\ot A)$ with this property.
\snl
In Section 3, we have the assumptions that
$$E(A\ot A)=\Delta(A)(1\ot A) \qquad\qquad\text{and}\qquad\qquad
         (A\ot A)E=(A\ot 1)\Delta(A),
$$
but as we  have $A^2=A$, then the assumption in A.1 will also be fulfilled. So, the conditions in A.1 are weaker than the ones used in Section 3 (as $A^2=A$ is assumed).
\nl
{\it In what follows}, we assume that $(A,\Delta)$ satisfies the assumption in 3.1 and that $E$ is the unique idempotent in $M(A\ot A)$ with this property.
\snl
We will prove the following (cf.\ Proposition 3.8 in Section 3).

\ainspr{A.2} Proposition \rm
There is a unique homomorphism $\Delta_1:M(A)\to M(A\ot A)$ that extends $\Delta$ from $A$ to $M(A)$ and that satisfies $\Delta_1(1)=E$.
\einspr

We will obtain the result from a more general result that we will prove first. This more general result will not only give the property in Proposition A.2 but it will also provide extensions of $\Delta\ot\iota$ and $\iota\ot\Delta$  needed later.

\ainspr{A.3} Proposition \rm
Let $A$ and $B$ be non-degenerate algebras and $\gamma:A\to M(B)$ a homomorphism. Assume that there is an idempotent element $e\in M(B)$ such that
$$\gamma(A)B=eB \qquad\qquad\text{and}\qquad\qquad B\gamma(A)=Be.$$
Then there is a unique homomorphism $\gamma_1:M(A)\to M(B)$, extending $\gamma$ and such that $\gamma_1(1)=e$.
\einspr

Before we prove this proposition, let us show the following.

\ainspr{A.4} Lemma \rm
With the assumptions of Proposition A.3 we have
$$e\gamma(a)=\gamma(a) \qquad\qquad\text{and}\qquad\qquad \gamma(a)e=\gamma(a)$$
for all $a\in A$.

\snl \bf Proof\rm: Take $a\in A$. For all $b\in B$ we have $\gamma(a)b\in eB$ and so $e\gamma(a)b=\gamma(a)b$. Because this is true for all $b\in B$, we have $e\gamma(a)=\gamma(a)$. Similarly $\gamma(a)e=\gamma(a)$.
\hfill$\square$
\einspr

As in the proof of Proposition 3.6, also here we will have that $e$ is the smallest idempotent with this property.
\snl
We now give {\it the proof of Proposition} A.3:

\ainspr{}Proof\rm: 
Assume first that $\gamma_1$ is such an extension. Take $m\in M(A)$ and $x,y\in B$. Then 
$$\gamma_1(m)x=\gamma_1(m)\gamma_1(1)x=\gamma_1(m)ex.$$
By assumption we can write $ex$ as a finite sum $\sum_i \gamma(a_i)b_i$ with $a_i\in A$ and $b_i\in B$ for all $i$ and then we get, using that $\gamma_1$ extends $\gamma$, 
$$\gamma_1(m)x=\sum_i \gamma(ma_i)b_i. \tag"(A.1)"$$
Similarly, 
$$y\gamma_1(m)=\sum_j c_j\gamma(d_jm) \tag"(A.2)"$$
if $ye$ is written as a finite sum $\sum_j c_j\gamma(d_j)$ with $c_j\in B$ and $d_j\in A$.
\snl
This already proves that such an extension is unique if it exists. It also suggests how to define this extension. To show that the above formulas (A.1) and (A.2) really define a multiplier $\gamma_1(m)$ for any multiplier $m\in M(A)$ one first has to argue e.g.\ that  
$$\sum_i \gamma(ma_i)b_i=0 \qquad\qquad\text{if}\qquad\qquad \sum_i \gamma(a_i)b_i=0$$
for a finite number of elements $a_i\in A$ and $b_i\in B$. 
\snl
To show this, take any $c\in B$ and $d\in A$. Then
$$c\gamma(d)\sum_i \gamma(ma_i)b_i = \sum_i c\gamma(dma_i)b_i=c\gamma(dm) \sum_i \gamma(a_i)b_i$$
and this is $0$ when we assume that $\sum_i \gamma(a_i)b_i=0$. By assumption $B\gamma(A)=Be$ and so it follows that also 
$$ye\sum_i \gamma(ma_i)b_i=0$$
for all $y\in B$. As $e\gamma(a)=\gamma(a)$ for all $a\in A$, we get $\sum_i \gamma(ma_i)b_i$ from the non-degeneracy of the product in $B$.
\snl
This shows that we can define $\gamma_1(m)x$ by the formula (A.1). Similarly, we can define $y\gamma_1(m)$ using the formula (A.2). Moreover, an argument as above will give that $y(\gamma_1(m)x)=(y\gamma_1(m))x$ so that indeed, $\gamma_1(m)$ is well-defined as an element in $M(B)$ for all $m\in M(A)$. 
\snl
Finally, it is not hard to show that $\gamma_1$ is still a homomorphism, that it extends $\gamma$ and that $\gamma_1(1)=e$.
This completes the proof.
\hfill$\square$\einspr

If $A$ and $B$ are $^*$-algebras and $\gamma$ is a $^*$-homomorphism, then the extension will still be a $^*$-homomorphism provided we have that $e^*=e$ in $M(B)$. 
But this property is a consequence of Lemma A.4, by taking adjoints and using that $e$ is the smallest idempotent with this property.
\snl
If $e=1$ we recover the original result about extending non-degenerate homomorphisms (as first shown in Proposition A.5 of [VD1]).
\snl
If we replace $B$ by $A\ot A$ and if we consider $\Delta$ for $\gamma$, then with  $e$ replaced by $E$, we get the result in Proposition A.2. 
\snl
Most of the time, we will use {\it the same symbol} for the extended homomorphisms. So, in particular, we will use $\Delta$ also for the extension $\Delta_1$ obtained in Proposition A.2.
\snl
If we replace $A$ by $A\ot A$ and $B$ by $A\ot A\ot A$ and if we consider the maps $\Delta\ot\iota$ and $\iota\ot\Delta$, we arrive at the following result.

\ainspr{A.5} Proposition \rm
The homomorphisms $\Delta\ot\iota$ and $\iota\ot\Delta$ have unique extensions to homomorphisms from $M(A\ot A)$ to $M(A\ot A\ot A)$, still denoted by $\Delta\ot\iota$ and $\iota\ot\Delta$ respectively, provided we require that
$$(\Delta\ot\iota)(1)=E\ot 1\qquad\qquad\text{and} \qquad\qquad (\iota\ot\Delta)(1)=1\ot E$$
where we use $1$ both for the identity in $M(A)$ and for the one in $M(A\ot A)$.
\einspr

Again this result is a straightforward application of the result in Proposition A.3. We just have to argue that 
$$(\Delta(A)\ot A)(A\ot A\ot A)=(E\ot 1)(A\ot A\ot A)$$
and similarly on the other side and for $\iota\ot\Delta$. Of course all of this follows from the assumptions $\Delta(A)(A\ot A)=E(A\ot A)$ and $(A\ot A)\Delta(A)=(A\ot A)E$ together with $A^2=A$ and similarly for the other one.
\nl
Let us now look at coassociativity. Using the extended maps $\Delta\ot \iota$ and $\iota\ot\Delta$ from Proposition A.5 we get the following expected formulation of coassociativity. 

\ainspr{A.6} Proposition \rm We have
$$(\Delta\ot\iota)\Delta(a)=(\iota\ot\Delta)\Delta(a)$$
for all $a\in A$.

\snl\bf Proof\rm: 
For clarity, we will use in the proof $\widetilde\Delta\ot \iota$ and $\iota\ot\widetilde\Delta$ for the extensions of $\Delta\ot\iota$ and $\iota\ot\Delta$ to the multiplier algebra $M(A\ot A)$.
\snl
Let $a\in A$.
We only need to argue that 
$$((\widetilde\Delta\ot\iota)\Delta(a))(1\ot 1\ot b)=(\Delta\ot\iota)(\Delta(a)(1\ot b))$$
for all $b\in A$. In a similar way we will get that 
$$(c\ot 1\ot 1)((\iota\ot\widetilde\Delta)\Delta(a))=(\iota\ot\Delta)((c\ot 1)\Delta(a))$$
for all $c \in A$ and then the result will follow from the formulation of coassociativity as in Definition 1.1.
\snl
To prove the claim, let $b\in B$ and $z\in A\ot A\ot A$. Write 
$$z(E\ot 1)=\sum_i u_i (\Delta\ot\iota)(v_i)$$
with $u_i\in A\ot A\ot A$ and $v_i\in A\ot A$. Then we get in a straightforward manner 
$$\align z((\widetilde\Delta\ot\iota)\Delta(a))(1\ot 1\ot b)
	&=z(E\ot 1)((\widetilde\Delta\ot\iota)\Delta(a))(1\ot 1\ot b)\\
	&=\sum_i u_i (\Delta\ot\iota)(v_i)((\widetilde\Delta\ot\iota)\Delta(a))(1\ot 1\ot b)\\
	&=\sum_i u_i ((\Delta\ot\iota)(v_i\Delta(a))(1\ot 1\ot b)\\
	&=\sum_i u_i (\Delta\ot\iota)(v_i\Delta(a)(1\ot b))\\
	&=\sum_i u_i (\Delta\ot\iota)(v_i)(\Delta\ot\iota)(\Delta(a)(1\ot b))\\
	&=z(E\ot 1)(\Delta\ot\iota)(\Delta(a)(1\ot b))\\
	&=z(\Delta\ot\iota)(\Delta(a)(1\ot b)).
\endalign$$
and this proves the claim. 
\hfill$\square$
\einspr
	
We now will push this result further and show that we also have the following (expected) formula.

\ainspr{A.7} Proposition \rm We have
$$(\Delta\ot\iota)\Delta(m)=(\iota\ot\Delta)\Delta(m)$$
for all $m\in M(A)$.
\einspr

Remark that from this result it would follow that $(\Delta\ot\iota)(E)=(\iota\ot \Delta)(E)$ because $E=\Delta(1)$. On the other hand, we want to prove this proposition using the uniqueness of the extensions of $(\Delta\ot\iota)\Delta$ and $(\iota\ot\Delta)\Delta$ on $A$. The result then will follow from the previous result. However, in order to obtain this uniqueness, we need the equality for $m=1$, that is, we need that $(\Delta\ot\iota)(E)=(\iota\ot \Delta)(E)$. 
\snl
Therefore, we {\it prove this equation first}.

\ainspr{A.8} Proposition \rm
We have
$(\Delta\ot\iota)(E)=(\iota\ot \Delta)(E)$.

\snl\bf Proof\rm:
Denote in this proof $(\Delta\ot\iota)(E)$ by $G$ and $(\iota\ot \Delta)(E)$ by $G'$.
We have
$$\align G(A\ot A\ot A)
	&= ((\Delta\ot\iota)(E))(A\ot A\ot A)\\
	&= ((\Delta\ot\iota)(E))(E\ot 1)(A\ot A\ot A)\\
   &= ((\Delta\ot\iota)(E))(\Delta(A)(A\ot A)\ot A)\\
	&= ((\Delta\ot\iota)(E(A\ot A)))(A\ot A\ot A)\\
	&= ((\Delta\ot\iota)\Delta(A))(A\ot A\ot A).
\endalign$$
Similarly we get 
$$G'(A\ot A\ot A)=((\iota\ot\Delta)\Delta(A))(A\ot A\ot A)$$
and because we know from Proposition A.6 that 
$$(\Delta\ot\iota)\Delta(A)=(\iota\ot\Delta)\Delta(A),$$
we find $G(A\ot A\ot A)=G'(A\ot A\ot A)$. In a similarly way we can show that 
$(A\ot A\ot A)G=(A\ot A\ot A)G'$ and as before, this implies that $G=G'$.
\hfill$\square$
\einspr

Now, the result in Proposition A.7 follows from the uniqueness in Proposition A.3.
\snl
We finish the discussion with an important remark.

\ainspr{A.9} Remark \rm
i) Because by definition we have $(\Delta\ot\iota)(1)=E\ot 1$, we also have
$$(\Delta\ot \iota)(E)=(E\ot 1)(\Delta\ot \iota)(E)= (\Delta\ot \iota)(E)(E\ot 1).$$
This means that 
$$(\Delta\ot \iota)(E)\leq E\ot 1.$$
Similarly 
$$(\iota\ot\Delta)(E)=(1\ot E)(\iota\ot\Delta)(E)=(\iota\ot\Delta)(E)(1\ot E)$$
and this means
$$(\iota\ot\Delta)(E)\leq 1\ot E.$$
Because the left hand sides are the same, we get an idempotent that is  smaller than both $E\ot 1$  and $1\ot E$.
\snl
ii) In the theory developed in this paper, we have that the two idempotents $E\ot 1$ and $1\ot E$ commute and that 
$$(\iota\ot\Delta)(E)=(1\ot E)(E\ot 1).$$
This means that the product of $E\ot 1$ and $1\ot E$ is again an idempotent in $M(A\ot A\ot A)$ and that $(\iota\ot\Delta)(E)$ is the biggest idempotent in $M(A\ot A\ot A)$ that is smaller than both $E\ot 1$ and $1\ot E$.
\hfill $\square$ \einspr

Let us now have a brief look at {\it examples}.

\aiinspr{A.10} Examples \rm
As we have already indicated, the assumption in A.1 is satisfied in all the examples we have considered up to now. This is shown in Example 3.6 for the two groupoid examples and the remark following 3.6 for the case of a weak Hopf algebra.
\snl
There is of course very little to say more in the case of a weak Hopf algebra as we start with an algebra with identity. In the two other cases, one can easily see what the extended maps are and that they satisfy the results formulated in this appendix. 
\snl
Finally, observe that in all these cases, the two idempotents $E\ot 1$ and $1\ot E$ commute and that 
$$(\iota\ot\Delta)(E)=(1\ot E)(E\ot 1).$$
This has been argued also in Section 3.
\hfill $\square$
\einspr

We have no examples where $E\ot 1$ and $1\ot E$ do not commute. We also have no examples where $(\Delta\ot \iota)(E)$ is not equal to $(E\ot 1)(1\ot E)$ and/or $(1\ot E)(E\ot 1)$. On the other hand, there seems to be no reason why these properties  should automatically hold.
\nl\nl



\bf References \rm
\nl
{[\bf A]} E.\ Abe: {\it Hopf algebras}. \rm Cambridge University Press (1977).
\snl
{[\bf B-N-S]} G.\ B\"ohm, F.\ Nill \& K.\ Szlach\'anyi: {\it Weak Hopf algebras I. Integral theory and C$^*$-structure}. J.\ Algebra 221 (1999), 385-438. 
\snl
{[\bf B-S]} G.\ B\"ohm  \& K.\ Szlach\'anyi: {\it Weak Hopf algebras II. Representation theory, dimensions and the Markov trace}. J.\ Algebra 233 (2000), 156-212. 
\snl
{[\bf Br]} R.\ Brown: {\it From groups to groupoids: A brief survey}. Bull. London Math. Soc. 19 (1987), 113-134.
\snl
{[\bf G]} K.R.\ Goodearl: {\it Von Neumann Regular Rings}. Pitman, London, 1979.
\snl
{[\bf H]} P.\ J.\ Higgins: {\it Notes on categories and groupoids}. Van Nostrand Reinhold, London (1971).
\snl
{[\bf N]} D.\ Nikshych: {\it On the structure of weak Hopf algebras}. Adv.\ Math.\ 170 (2002), 257-286.
\snl
{[\bf N-V1]} D.\ Nikshych \& L.\ Vainerman: {\it Algebraic versions of a finite dimensional quantum groupoid}. Lecture Notes in Pure and Applied Mathematics 209 (2000), 189-221. 
\snl
{[\bf N-V2]} D.\ Nikshych \& L.\ Vainerman: {\it Finite quantum groupiods and their applications}. In {\it New Directions in Hopf algebras}. MSRI Publications, Vol.\ 43 (2002), 211-262.
\snl
{[\bf P]} A.\ Paterson: {\it Groupoids, inverse semi-groups and their operator algebras}. Birkhauser, Boston (1999).
\snl
{[\bf R]} J.\ Renault: {\it A groupoid approach to C$^*$-algebras}. Lecture Notes in Mathematics 793, Springer Verlag.
\snl
{[\bf S]} M.\ Sweedler: {\it Hopf algebras}. Benjamin, New-York (1969).
\snl
{[\bf V]} L.\ Vainerman (editor): {\it Locally compact quantum groups and groupoids}. IRMA Lectures in Mathematics and Theoretical Physics 2, Proceedings of a meeting in Strasbourg, de Gruyter (2002).
\snl
{[\bf VD1]} A.\ Van Daele: {\it Multiplier Hopf algebras}. Trans. Am. Math. Soc.  342(2) (1994), 917-932.
\snl
{[\bf VD2]} A.\ Van Daele: {\it An algebraic framework for group duality}. Adv. in Math. 140 (1998), 323-366.
\snl
{[\bf VD3]} A.\ Van Daele: {\it Tools for working with multiplier Hopf algebras}. Preprint University of Leuven. Arxiv math.QA/0806.2089. ASJE (The Arabian Journal for Science and Engineering) C - Theme-Issue   33 (2008), 505--528. 
\snl
{[\bf VD4]} A.\ Van Daele: {\it Separability and multiplier algebras}. Preprint University of Leuven (2012 - in preparation).
\snl 
{[\bf VD-Ve]} A.\ Van Daele \& J.\ Vercruysse: {\it Local units and multiplier algebras}. University of Leuven and University of Brussels (2012 - in preparation).  
\snl
{[\bf VD-W1]} A.\ Van Daele \& S.\ Wang: {\it Multiplier Unifying Hopf Algebras}. Preprint University of Leuven and Southeast University of Nanjing (2008 - unpublished).
\snl 
{[\bf VD-W2]} A.\ Van Daele \& S.\ Wang: {\it Weak multiplier Hopf algebras I. The main theory}. Preprint University of Leuven and Southeast University of Nanjing (October 2012).
\snl 
{[\bf VD-W3]} A.\ Van Daele \& S.\ Wang: {\it Weak multiplier Hopf algebras II. The source and target algebras}. Preprint University of Leuven and Southeast University of Nanjing (2012).
\snl
{[\bf VD-W4]} A.\ Van Daele \& S.\ Wang: {\it Weak multiplier Hopf algebras III. Integrals and duality}. Preprint K.U.\ Leuven and Southeast University of Nanjing (2012 - in preparation).
\snl 
{[\bf Ve]} J.\ Vercruysse: {\it Local units versus local projectivity dualisations: Corings with local structure maps}. Commun. in Alg. 34 (2006) 2079-2103.
\snl

\end